\documentclass{amsart}
\usepackage{amsfonts}

\usepackage{graphicx}
\usepackage{amscd}

\newtheorem{theorem}{Theorem}
\theoremstyle{plain}
\newtheorem{acknowledgement}{Acknowledgement}

\newtheorem{definition}{Definition}

\newtheorem{proposition}{Proposition}
\newtheorem{remark}{Remark}

\numberwithin{equation}{section}
\newcommand{\thmref}[1]{Theorem~\ref{#1}}
\newcommand{\secref}[1]{\S\ref{#1}}
\newcommand{\lemref}[1]{Lemma~\ref{#1}}
\input{tcilatex}

\begin{document}
\title[$q$-Bernoulli Numbers and Polynomials...]{$q$-Bernoulli Numbers and
Polynomials Associated with Multiple $q$-Zeta Functions and Basic $L$-series}
\author{TAEKYUN KIM}
\address{TAEKYUN KIM\\
Institute of Science Education, Kongju National University Kongju 314-701,
S. Korea}
\email{tkim@kongju.ac.kr}
\author{YILMAZ\ SIMSEK}
\curraddr{YILMAZ\ SIMSEK\\
Mersin University, Faculty of Science, Department of Mathematics 33343
Mersin, Turkey}
\email{ysimsek@mersin.edu.tr}
\author{H. M. SRIVASTAVA}
\address{H. M. SRIVASTAVA\\
Department of Mathematics and statistics University of Victoria, Victoria
British Columbia V8W 3P4 Canada}
\email{harimsri@math.uvic.ca}

\begin{abstract}
By using $q$-Volkenborn integration and uniform differentiable on $\mathbb{Z}%
_{p}$, we construct $p$-adic $q$-zeta functions. These functions interpolate
the $q$-Bernoulli numbers and polynomials. The value of $p$-adic $q$-zeta
functions at negative integers are given explicitly. We also define new
generating functions of $q$-Bernoulli numbers and polynomials. By using
these functions, we prove analytic continuation of some basic (or $q$- ) $L$%
-series. These generating functions also interpolate Barnes' type Changhee $%
q $-Bernoulli numbers with attached to Dirichlet character as well. By
applying Mellin transformation, we obtain relations between Barnes' type $q$%
-zeta function and new Barnes' type Changhee $q$-Bernolli numbers.
Furthermore, we construct the Dirichlet type Changhee ( or $q$-) $L$%
-functions.
\end{abstract}

\maketitle

\section{Introduction, Definition and Notations}

For any complex number $z$, it is well known that the usual Bernoulli
polynomials $B_{n}(z)$ are defined by means of of the generating function
(see \cite{E. T. Wittaker and G. N. Watson}, \cite{K. Shiratani and S.
Yamamoto}, \cite{Apostol}, and \cite{Deeba}):
\begin{equation}
F(t,x)=\frac{te^{zt}}{e^{t}-1}=\sum_{n=0}^{\infty }B_{n}(z)\frac{t^{n}}{n!}%
\text{ \ (}\mid t\mid <2\pi \text{ ).}  \label{eq-1}
\end{equation}%
Note that, substituting $z=0$ into (\ref{eq-1}), $B_{n}(0)=B_{n}$ is the
usual $n$th Bernoulli number:
\begin{equation}
F(t)=\frac{t}{e^{t}-1}=\sum_{n=0}^{\infty }B_{n}\frac{t^{n}}{n!}\text{ \ (}%
\mid t\mid <2\pi \text{ ).}  \label{q-2}
\end{equation}

Over five decades ago, Carlitz \cite{Carlitz} defined $q$-extensions of
these classical Bernoulli numbers and polynomials and proved properties
generating those satisfied by $B_{n}$ and $B_{n}(z).$ Recently, Koblitz \cite%
{N. Koblitz} used these properties, especially the so-called \textit{%
distribution relation} for $q$-Bernoulli polynomials, in order to construct
the corresponding $q$-extensions of the $p$-adic measures and to define a $q$%
-extension of $p$-adic Dirichlet $L$-series.

When one talks of $q$-extensions, $q$ can be variously considered as an
indetermined, a complex number $q\in \mathbb{C}$, or, when $p$ be a prime
number, a $p$-adic number $q\in \mathbb{C}_{p}$, where $\mathbb{C}_{p}$ is
the $p$-adic completion of the algebraic closure of $\mathbb{Q}_{p}$. If $%
q\in \mathbb{C}$, one normally assumes $\mid q\mid <1$. If $q\in \mathbb{C}%
_{p}$, then we assume
\begin{equation*}
\mid q-1\mid _{p}<p^{-\frac{1}{p-1}},
\end{equation*}%
so that%
\begin{equation*}
q^{x}=\exp (x\log q)\text{ for}\mid x\mid _{p}\leq 1.
\end{equation*}%
We use the notation (see also \cite{H. M. Srivastava and P. W. Karlsson}, p.
346 \textit{et seq.}):
\begin{equation*}
\lbrack x]=[x:q]=\frac{1-q^{x}}{1-q}.
\end{equation*}%
Thus%
\begin{equation*}
\lim_{q\rightarrow 1}[x:q]=x
\end{equation*}%
for any $x\in \mathbb{C}$ in the complex case and any $x$ with $\mid x\mid
_{p}\leq 1$ in the $p$-adic case (see \cite{N. Koblitz}, \cite{N. Koblitz-1}
and \cite{Kim4}).

Carlitz's $q$-Bernoulli numbers $\beta _{n}=\beta _{n}(q)$ can be determined
inductively by\cite{Carlitz}
\begin{equation*}
\beta _{0}=1,q(q\beta +1)^{n}-\beta _{n}=\left\{
\begin{array}{c}
1,\text{ if }n=1 \\
0,\text{ if }n>1,%
\end{array}
\right.
\end{equation*}
with the usual convention about replacing $\beta ^{n}$ by $\beta _{n}$.

The $q$-Bernoulli polynomials $\beta _{n}(x:q)$ are given as $(q^{x}\beta
+[x])^{n}$, i.e., as follows
\begin{equation*}
\beta _{n}(x:q)=\sum_{k=0}^{n}\left(
\begin{array}{c}
n \\
k%
\end{array}%
\right) \beta _{k}q^{kx}[x]^{n-k}.
\end{equation*}%
As $q\rightarrow 1$, we have $\beta _{n}(q)\rightarrow B_{n},\beta
_{n}(x:q)\rightarrow B_{n}(x).$

Let $\chi $ be a Dirichlet character of conductor $f\in \mathbb{Z}^{+}$, the
set of positive integer numbers. Then the generalized Bernoulli numbers $%
B_{n,\chi }$ are defined by
\begin{equation}
F_{\chi }(t)=\sum_{a=0}^{f-1}\frac{\chi (a)te^{at}}{e^{ft}-1}%
=\sum_{n=0}^{\infty }B_{n,\chi }\frac{t^{n}}{n!}\text{ \ (}\mid t\mid <2\pi
\text{ ).}  \label{eq-3}
\end{equation}%
Then we defined generalized Carlitz's $q$-Bernoulli number $\beta _{m,\chi
}(q)=\beta _{m,\chi }$ as follows\cite{Iwasawa}, \cite{Kim1}, \cite{Howard}
\begin{equation}
\beta _{m,\chi }(q)=[f]^{m-1}\sum_{a=0}^{f-1}\chi (a)q^{a}\beta _{m}(\frac{a%
}{f}:q^{f}).  \label{eq-4}
\end{equation}%
As $q\rightarrow 1$,(\ref{eq-4}) is reduced to (\ref{eq-3}).

The Euler numbers $E_{n}$ are usually defined by means of of the following
generating function (see, for example, \cite{H. M. Srivastava and J. Choi},
p. 63, Eq. 1.6 (40); see also \cite{H. Tsumura-1}) for different
definition):
\begin{equation*}
\frac{2e^{t}}{e^{2t}+1}=\sec h(t)=\sum_{n=0}^{\infty }E_{n}\frac{t^{n}}{n!}%
\text{ \ (}\mid t\mid <\frac{\pi }{2}\text{ ).}
\end{equation*}%
These numbers are classical and important in number theory. Frobenius
extended such numbers as $E_{n}$ to the so-called Frobenius-Euler numbers $%
H_{n}(u)$ belonging to an algebraic number $u$, with $\mid u\mid >1$, and
many authors have investigated their properties ( \cite{Kim4}, \cite{Kim9}
). Shiratani and Yamamoto \cite{K. Shiratani and S. Yamamoto} constructed a $%
p$-adic interpolation $G_{p}(s,u)$ of the Frobenius-Euler numbers $H_{n}(u)$
and as its application, they obtained an explicit formula for $%
L_{p}^{^{\prime }}(0,\chi )$ with any Dirichlet character $\chi $. In \cite%
{H. Tsumura-1}, Tsumura defined the generalized Frobenius-Euler numbers $%
H_{n,\chi }(u)$ for any Dirichlet character $\chi $, which are analogous to
the generalized Bernoulli numbers. He constructed their\ Shiratani and
Yamamoto $p$-adic interpolation $G_{p}(s,u)$ of $H_{n}(u).$

Let $u$ be an algebraic number. For $u\in \mathbb{C}$ with $|u|>1$, the
Frobenius-Euler numbers $H_{n}(u)$ belonging to $u$ are defined by means of
of the generating function
\begin{equation*}
\frac{1-u}{e^{t}-u}=e^{H(u)t}
\end{equation*}%
with usual convention of symbolically replacing $H^{n}(u)$ by $H_{n}(u)$.
Thus for the Frobenius-Euler numbers $H_{n}(u)$ belonging to $u$, we have (
see\cite{K. Shiratani})
\begin{equation}
\frac{1-u}{e^{t}-u}=\sum_{n=0}^{\infty }H_{n}(u)\frac{t^{n}}{n!}.
\label{eq-5}
\end{equation}%
By using (\ref{eq-5}), and following the usual convention of symbolically
replacing $H^{n}(u)$ by $H_{n}(u)$, we have%
\begin{equation*}
\text{ }H_{0}=1\text{ and }(H(u)+1)^{n}=uH_{n}(u)\text{ for (}n\geq 1\text{).%
}
\end{equation*}%
We also note that
\begin{equation*}
H_{n}(-1)=\mathfrak{E}_{n},
\end{equation*}%
where $\mathfrak{E}_{n}$\ denotes the aforementioned Tsumura version ( see%
\cite{K. Shiratani}) of the classical Euler numbers $E_{n}$ which we
recalled above.

\bigskip For an algebraic number $u\in \mathbb{C}$ with $|u|>1$, the
Frobenius-Euler polynomials belonging to $u$, that is, the polynomials $%
H_{n}(u,x)$ are defined by (see\cite{H. Tsumura-1})
\begin{equation}
\frac{1-u}{e^{t}-u}e^{xt}=e^{H(u,x)t}=\sum_{n=0}^{\infty }H_{n}(u,x)\frac{%
t^{n}}{n!}  \label{eq-6}
\end{equation}%
with usual convention of symbolically replacing $H^{n}$ by $H_{n}$ as
before. By using (\ref{eq-5}) and (\ref{eq-6}), we readily have
\begin{equation*}
H_{n}(u,0)=H_{n}(u)\text{ and }H_{n}(u,x)=\sum_{k=0}^{n}\left(
\begin{array}{c}
n \\
k%
\end{array}%
\right) H_{k}(u)x^{n-k}.
\end{equation*}

Let $\chi $ be a Dirichlet character of conductor $f\in \mathbb{Z}^{+}$. We
define the $n$th generalized Euler numbers $H_{n,\chi }(u)$ belonging to $u$%
, by\cite{H. Tsumura-1}
\begin{equation}
\sum_{a=0}^{f-1}\frac{(1-u^{f})\chi (a)e^{at}u^{f-a-1}}{e^{ft}-u^{f}}%
=\sum_{n=0}^{\infty }H_{n,\chi }(u)\frac{t^{n}}{n!}.  \label{eq-7}
\end{equation}%
By using (\ref{eq-5}) to (\ref{eq-7}), we can easily see that
\begin{eqnarray*}
H_{n,\chi }(u) &=&f^{n}\sum_{a=0}^{f-1}\chi (a)u^{f-a-1}H_{n}(u^{f},\frac{a}{%
f}) \\
&=&\sum_{a=0}^{f-1}\chi (a)u^{f-a-1}\sum_{k=0}^{n}\left(
\begin{array}{c}
n \\
k%
\end{array}%
\right) H_{k}(u^{f})a^{n-k}f^{k}.
\end{eqnarray*}%
We note that, when $\chi =1$, we have
\begin{equation*}
H_{n,1}(u)=H_{n}(u)\text{, for ( }n\geq 0\text{ ).}
\end{equation*}

Carlitz\cite{Carlitz} also defined $q$-Euler numbers and polynomials as
follows:
\begin{equation*}
H_{0}(u:q)=1\text{ and }(qH(u:q)+1)^{n}-uH_{n}(u:q)=0\text{ ( }n\geq 1\text{%
),}
\end{equation*}%
where $u$ is a complex number $\mid u\mid >1$. For $n\geq 0$,\ and with the
usual convention of replacing $H^{n}$ by $H_{n}$,\ we have (see \cite%
{Carlitz})%
\begin{equation}
H_{n}(u,x:q)=(q^{x}H(u,x:q)+[x])^{n}.  \label{eq-8}
\end{equation}

When $q\rightarrow 1$ in (\ref{eq-8}), we obtain the following limit
relationship with the Frobenius-Euler numbers $H_{n}(u)$ (\ref{eq-5}):%
\begin{equation*}
\lim_{q\rightarrow 1}H_{n}(u,1:q)=H_{n}(u)
\end{equation*}%
( see, for detail \cite{Iwasawa}, \cite{Dilcher}, \cite{Andrews},\cite{Kim9}%
).

Consider the finite products $E_{n}=\prod_{j<n}X_{j\text{ }}$ of a sequence $%
(X_{j})_{j\geq 0}$ of sets. We would like to say that these partial products
converge to the infinite product $E=\prod_{j\geq 0}X_{j\text{ }}$\ and thus
consider this last product as limit of the sequence $(E_{n})$. The \textit{%
projective limit} $E=\lim_{\leftarrow }E_{n}$ is defined by (see \cite{A. M.
Robert}, p. 28):

\begin{definition}
A sequence $(E_{n},\varphi _{n})_{n\geq 0}$ of sets and maps $\varphi
_{n}:E_{n+1}\rightarrow E_{n}$ ($n\geq 0$) is called a projective system. A
set $E$ given to gether with maps $\psi _{n}:E\rightarrow E_{n}$ such that $%
\psi _{n}=\varphi _{n}o\psi _{n+1}$ ($n\geq 0$) is called a projective limit
of the sequence $(E_{n},\varphi _{n})_{n\geq 0}$ if the following condition
is satisfied:

For each set $X$ and maps $f_{n}:X\rightarrow E_{n}$ satisfying $%
f_{n}=\varphi _{n}of_{n+1}$ ($n\geq 0$) there is a unique factorization $f$
of $f_{n\text{ }}$through the set $E$:%
\begin{equation*}
f_{n}=\psi _{n}of:X\rightarrow E\rightarrow E_{n\text{ }}\text{ }(n\geq 0).
\end{equation*}%
The maps $\varphi _{n}:E_{n+1}\rightarrow E_{n}$ are usually called
transition maps of the projective system. The whole system, represented by%
\begin{equation*}
E_{0}\leftarrow E_{1}\leftarrow ...\leftarrow E_{n}\leftarrow ...,
\end{equation*}%
is also called an inverse system.
\end{definition}

Kim \cite{Kim5} defined the $q$-Volkenborn integration and gave relations
between the $q$-Bernoulli numbers and the $q$-Euler numbers. He constructed
a new measure as well.

Let $p$ be a fixed prime. For a fixed positive integer $d$ with $(p,d)=1$,
we set (see \cite{Kim5})

\begin{eqnarray*}
X &=&X_{d}=\lim_{\leftarrow _{N}}\mathbb{Z}/\mathbb{Z}dp^{N}, \\
\text{ }X_{1} &=&\mathbb{Z}_{p}, \\
\text{ }X^{\ast } &=&\cup _{%
\begin{array}{c}
0<a<dp \\
(a,p)=1%
\end{array}%
}a+dp\mathbb{Z}_{p}
\end{eqnarray*}%
and%
\begin{equation*}
a+dp^{N}\mathbb{Z}_{p}=\left\{ x\in X\mid x\equiv a(\func{mod}%
dp^{N})\right\} ,
\end{equation*}%
where $a\in \mathbb{Z}$ satisfies the condition $0\leq a<dp^{N}$(\cite{Kim3}%
). The $p$-adic absolute value in is normalized in such a way that
\begin{equation*}
\mid p\mid _{p}=\frac{1}{p}.
\end{equation*}

We say that $f$ is a uniformly differentiable function at a point $a\in
\mathbb{Z}_{p}$, and write $f\in UD(\mathbb{Z}_{p})$, if the difference
quotient
\begin{equation*}
F_{f}(x,y)=\frac{f(x)-f(y)}{x-y}
\end{equation*}%
has a limit
\begin{equation*}
l=f^{^{\prime }}(a)\text{ as }(x,y)\rightarrow (a,a).
\end{equation*}%
For $f\in UD(\mathbb{Z}_{p})$, let us begin with the expression:
\begin{equation*}
\frac{1}{[p^{N}]}\sum_{0\leq j<p^{N}}q^{j}f(j)=\sum_{0\leq j<p^{N}}f(j)\mu
_{q}(j+p^{N}\mathbb{Z}_{p}),
\end{equation*}%
which represents a $q$analogue of Riemann sums for $f$. The integral of $f$
on $\mathbb{Z}_{p}$ is defined as the limit of these sums ( as $N\rightarrow
\infty $ ) if this limit exists. The $q$-Volkenborn integral of a function $%
f\in UD(\mathbb{Z}_{p})$ is defined by
\begin{equation}
\int_{\mathbb{Z}_{p}}f(x)d\mu _{q}(x)=\lim_{N\rightarrow \infty }\frac{1}{%
[p^{N}]}\sum_{0\leq j<p^{N}}q^{j}f(j),  \label{Eq-17}
\end{equation}%
where
\begin{eqnarray*}
\mu _{q}(j) &=&\mu _{q}(j+p^{N}\mathbb{Z}_{p})=\frac{q^{j}}{[p^{N}]} \\
(0 &\leq &j<p^{N};\text{ }N\in \mathbb{Z}^{+})
\end{eqnarray*}
(see, for detail \cite{Kim3}, \cite{Kim4}, \cite{Kim2}, \cite{Kim6}, \cite%
{Kim7}).

Kurt Hensel (1861-1941) invented the so-called $p$-adic numbers around the
end of the nineteenth century. In sipite of their being already one hundred
years old, these numbers are still today enveloped in an aura of mystery
within scientific community. Although they have penetrated several
mathematical fields, Number Theory, Algebraic Geometry, Algebraic Topology,
Analysis, Mathematical Physics, String Theory, Field Theory, Stochastic
Differential Equations on real Banach Spaces and Manifolds and other parts
of the natural sciences, they have yet to reveal their full potentials in
(for example) physics. While solving mathematical and physical problems and
while constructing and investigating measures on manifolds, the $p$-adic
numbers are used. There is an unexpected connection of the $p$-adic Analysis
with $q$-Analysis and Quantum Groups, and thus with Noncommutative Geometry,
and $q$-Analysis is a sort of $q$-deformation of the ordinary analysis.
Spherical functions on Quantum Groups are $q$-special functions. ( see \cite%
{Kim5}, \cite{Kim7}, \cite{Kim9}, \cite{Kim10}, \cite{Kim3}, \cite{A. M.
Robert}, \cite{Katriel}, \cite{W. H. Schikhof}, \cite{Vilademir}, \cite%
{Washington}, \cite{Khrennikov}, \cite{Khrennikov-1} ).

Kim\cite{Kim7} defined the Daehee numbers, $D_{m}(z:q)$ by using an
invariant $p$-adic $z$-integrals as follows:
\begin{equation*}
\mu _{z}(a+p^{N}\mathbb{Z}_{p})=\frac{z^{a}}{[p^{N}:z]}
\end{equation*}%
and%
\begin{equation*}
\lbrack x:z]=\frac{1-z^{x}}{1-z},
\end{equation*}%
which can be extended to distributions on $\mathbb{Z}_{p}$,
\begin{equation}
D_{m}(z:q)=\int_{\mathbb{Z}_{p}}[x]^{m}d\mu _{z}(x),  \label{Eq-10}
\end{equation}%
$z\in \mathbb{C}_{p}$. If we take $z=q$ in (\ref{Eq-10}), then we observe
that%
\begin{equation*}
D_{m}(q:q)=\beta _{m}(q)
\end{equation*}%
in terms of Carlitz's $q$-Bernoulli numbers mentioned above. In the case
when $z=u$ in (\ref{Eq-10}), the Daehee numbers become the $q$-Eulerian
numbers as follows:
\begin{equation*}
D_{m}(u:q)=\int_{\mathbb{Z}_{p}}[x]^{m}d\mu _{u}(x)=H_{m}(u^{-1}:q).
\end{equation*}%
By the definition of the Daehee numbers, we easily see that
\begin{equation*}
D_{m}(u:q)=\frac{1}{(1-q)^{m}}\sum_{l=0}^{m}\left(
\begin{array}{c}
m \\
l%
\end{array}%
\right) (-1)^{m-l}\frac{l+1}{[l+1]}.
\end{equation*}%
The Daehee polynomials are defined as follows\cite{Kim7}, \cite{Kim4}:
\begin{equation}
D_{m}(z,x:q)=\int_{\mathbb{Z}_{p}}[x+t]^{m}d\mu _{z}(t),  \label{Eq-11}
\end{equation}%
$z\in \mathbb{C}_{p}$.

\bigskip We readily see from (\ref{Eq-11}) that
\begin{equation*}
D_{m}(z,x:q)=\sum_{l=0}^{m}\left(
\begin{array}{c}
m \\
l%
\end{array}%
\right) q^{lx}[x]^{n-l}D_{l}(z:q)=(q^{x}D(z:q)+[x])^{m}
\end{equation*}%
with usual convention of symbolically replacing $D^{m}(z,x:q)$ by $%
D_{m}(z,x:q)$.

We note also that
\begin{equation*}
D_{m}(z,x:q)=H_{m}(u^{-1},x:q),
\end{equation*}%
and%
\begin{equation*}
D_{m}(q,x:q)=\beta _{m}(x:q),
\end{equation*}%
where $H_{m}(u^{-1},x:q)$ and $\beta _{m}(x:q)$ are Carlitz's $q$-Euler
Polynomials and Carlitz's $q$-Bernoulli Polynomials, respectively.

Ruijsenaars\cite{S. N. M. Ruijsenaars} showed how various known results
concerning the Barnes multiple zeta and gamma functions can be obtained as
specializations of simple features shared by a quite extensive class of
functions. The pertinent functions involve Laplace transforms, and their
asymptotic was obtained by exploiting this. He demonstrated how Barnes'
multiple zeta and gamma functions fit into a recently developed theory of
minimal solutions to first-order analytic difference equations. Both of
these approaches to the Barnes functions gave rise to novel integral
representations.

In one of an impressive series of papers ( \cite{Barnes}; see also \cite{H.
M. Srivastava and J. Choi}, Chapter 2), Barnes developed the so-called
multiple zeta and multiple gamma functions. Barnes' multiple zeta function $%
\zeta _{N}(s,w\mid a_{1},...,a_{N})$ depends on parameters $a_{1},...,a_{N}$
that will be taken positive throughout this paper. It defined by the series
\begin{equation}
\zeta _{N}(s,w\mid a_{1},...,a_{N})=\sum_{m_{1,...,m_{N}=0}}^{\infty
}(w+m_{1}a_{1}+...+m_{N}a_{N})^{-s},  \label{Eq-12}
\end{equation}%
where $\func{Re}(w)>0,\func{Re}(s)>N$. From the definition (\ref{Eq-12}), we
immediately obtain the recurrence relation \cite{S. N. M. Ruijsenaars}:
\begin{equation}
\zeta _{M+1}(s,w+a_{M+1}\mid a_{1},...,a_{N+1})-\zeta _{M+1}(s,w\mid
a_{1},...,a_{N+1})=-\zeta _{M}(s,w\mid a_{1},...,a_{N})  \label{eq-9}
\end{equation}%
with
\begin{equation*}
\zeta _{0}(s,w)=w^{-s}.
\end{equation*}

Barnes showed that $\zeta _{N}$ has a meromorphic continuation in $s$ ( with
simple poles only at $s=1,...,N$ ) and defined his multiple gamma function $%
\Gamma _{N}^{B}(w)$ in terms of the $s$-derivative at $s=0$, which may be
recalled here as follows\cite{S. N. M. Ruijsenaars}:
\begin{equation*}
\Psi _{N}(w\mid a_{1},...,a_{N})=\partial _{s}\zeta _{N}(s,w\mid
a_{1},...,a_{N})\mid _{s=0}.
\end{equation*}%
Clearly, analytic continuation of (\ref{eq-9}) yields the recurrence
relation:
\begin{eqnarray*}
\Psi _{M+1}(w+a_{M+1} &\mid &a_{1},...,a_{N+1})-\Psi _{M+1}(w\mid
a_{1},...,a_{N+1}) \\
&=&-\Psi _{M}(w\mid a_{1},...,a_{N}),
\end{eqnarray*}%
with%
\begin{equation*}
\Psi _{0}(w)=-\ln w.
\end{equation*}

Up to inessential factors, the functions\ $\zeta _{1}$ and $\Psi _{1}$ are
equal to the Hurwitz zeta function and the $\Psi $\ (or the digamma)
function (cf. e.g., Ref. \cite{E. T. Wittaker and G. N. Watson}). For $%
a_{1}=a_{2}=1$, the function
\begin{equation*}
\exp (\Psi _{2}(a_{1}+a_{2}-w\mid a_{1},a_{2})-\Psi _{2}(w\mid a_{1},a_{2}))
\end{equation*}%
was already studied by H\"{o}lder in 1886 \cite{S. N. M. Ruijsenaars}. It
was called the double sine function by Kurokawa\cite{N. Kurokawa}. In fact,
Kurokawa\cite{N. Kurokawa} considered multiple sine functions defined in
terms of $\Psi _{N}(w)$ and related these functions to Selberg zeta
functions and the determinants of Laplacians occruing in symmetric space
theory \cite{N. Kurokawa}. Barnes' multiple zeta and gamma functions were
also encountered by Shintani within the context \ of analytic number theory.
In recent years, they showed up in the form actor program for integrable
field theories and in studies of XXZ model correlation functions\cite{M.
Jimbo and T. Miwa}. See also recent paper by M. Nishizawa \cite{M. Nishizawa}%
, where $q$-analogues of the multiple gamma functions are studied. Friedman
and Ruijsenaars\cite{E. Friedman and S. Ruijsenaars} showed that Shintani's
work on multiple zeta and gamma functions can be simplified and extended by
exploiting difference equations. They re-proved many of Shintani's formulas
and prove several new ones. They also relate Barnes' triple gamma function
to the elliptic gamma function appearing in connection with certain
integrable systems.

The Barnes' multiple Bernoulli polynomials $B_{n}(x,r\mid a_{1},...,a_{r})$,
cf. \cite{Barnes}, are defined by
\begin{equation}
\frac{t^{r}e^{xt}}{\prod_{j=1}^{r}(e^{a_{j}t}-1)}=\sum_{n=0}^{\infty
}B_{n}(x,r\mid a_{1},...,a_{r})\frac{t^{n}}{n!},  \label{Eq-13}
\end{equation}
for $\mid t\mid <1$.

By (\ref{Eq-12}) and (\ref{Eq-13}), it is easy to see that
\begin{equation*}
\zeta _{N}(-m,w\mid a_{1},...,a_{N})=\frac{(-1)^{N}m!}{(N+m)!}%
B_{N+m}(w,N\mid a_{1},...,a_{N}),
\end{equation*}
for $w>0$ and $m$ is a positive integer.

In recent years, many mathematicians and physicians have investigated zeta
functions, multiple zeta functions, $L$-series, and multiple $q$-Bernoulli
numbers and polynomials because mainly of their interest and importance.
These functions and numbers are in used not only in Complex Analysis and
Mathematical Physics, but also in used in $p$-adic Analysis and other areas.
In particular, multiple zeta functions occur within the context of Knot
Theory, Quantum Field Theory, Applied Analysis and Number Theory (see \cite%
{Askey}, \cite{Barnes}, \cite{Cherednik}, \cite{Kim3}, \cite{Kim5}, \cite%
{Kim7}, \cite{Kim9},\cite{kim-Jang-Rim-Son}, \cite{T. M. Rassia and H. M.
Srivastava}, \cite{Y. Simsek}, \cite{C. A. Nelson and M. G. Gartley}, \cite%
{C. A. Nelson and M. G. Gartley-1}, \cite{P. T. Young}, \cite{B. E. Sagan}, %
\cite{N. Koblitz-1}, \cite{T. H. Koornwinder} ).

Kim\cite{Kim10} studied on the multiple $L$-series and functional equation
of this functions. He found the value of this function at negative integers
in terms of generalized Bernoulli numbers.

Russias and Srivastava\cite{T. M. Rassia and H. M. Srivastava} presented a
systematic investigation of several families of infinite series which are
associated with the Riemann zeta functions, the digamma functions,the
harmonic numbers, and the Stirling numbers of the first kind.

Matsumoto\cite{K. Matsumoto} considered general multiple zeta functions of
several variables, involving both Barnes multiple zeta functions and
Euler-Zagier sums as special cases. He proved the meromorphic continuation
to the whole space , asymptotic expansions, and upper bounded estimates.
These results were expected to have applications to some arithmetical $L$%
-functions. His method was based on \ the classical Mellin-Barnes integral
formula.

Ota\cite{K. Ota} studied on Kummer-type congruences for derivatives of
Barnes' multiple Bernoulli polynomials. Ota\cite{K. Ota} also generalized
these congruences to derivatives of Barnes' multiple Bernoulli polynomials
by an elementary method and gave a $p$-adic interpolation of them.
Subsequently, Ota \cite{K. Ota-1} defined derivatives of the Dedekind sums
and their reciprocity law. They were obtained from values at non-positive
integers of the first derivatives of Barnes' double zeta functions. As
special cases, they give finite product expressions of the Stirling modular
form and the double gamma function at positive rational numbers. The
original Dedekind sum appears at various places in mathematics, so the
derivative of the Dedekind sums may be expected to be useful as well. It
would be very interesting if we could obtain different proofs from Ota's for
the reciprocity laws about derivatives, just as the original reciprocity law
of Dedekind was obtained from the transformation formulas of $\log \eta (z)$%
. We note that by considering Barnes' $r$-ple zeta function \cite{K. Ota-1}
or zeta functions with characters Ota obtain reciprocity laws for sums
involving derivatives of the Barnes $r$-ple Bernoulli polynomials or
Dedekind sums with character.

Simsek\cite{Y. Simsek-1} gave relations between zeta functions,
trigonometric functions and Dedekind sums. He also found reciprocity law of
this sums related to Lambert series and Eisenstein series.

Woon\cite{S. C. Woon} presented a series of diagrams showing the Julia set
of the Riemann zeta functions and its related Mandelbrot set. The Julia and
Mandelbrot sets of the Riemann zeta function have unique features and are
quite unlike those of any elementary functions.

By using non-Archimedean $q$-integration, Kim \cite{Kim7} introduced
multiple Changhee $q$-Bernoulli polynomials which form a $q$-extension of
Barnes' multiple Bernoulli polynomials. He also constructed the Changhee $q$%
-zeta functions (which gives $q$-analogues of Barnes' multiple zeta
functions) and indicated some relationships between the Changhee $q$-zeta
function and Daehee $q$-zeta function.

A sequence of $p$-adic rational numbers as multiple Changhee $q$-Bernoulli
numbers and polynomials are defined as follows\cite{Kim7}, \cite{Kim14}:

Let $a_{1},...,a_{r}$ be nonzero elements of the $p$-adic number field and
let $z\in \mathbb{C}_{p}$.
\begin{equation}
\beta _{n}^{(r)}(w:q\mid a_{1},...,a_{r})=\frac{1}{\prod_{j=1}^{r}a_{j}}%
\int_{\mathbb{Z}_{p}^{r}}[w+\sum_{j=1}^{r}a_{j}x_{j}]^{n}d\mu _{q}(x),
\label{Eq-14}
\end{equation}%
and
\begin{equation*}
\beta _{n}^{(r)}(q\mid a_{1},...,a_{r})=\frac{1}{\prod_{j=1}^{r}a_{j}}\int_{%
\mathbb{Z}_{p}^{r}}[\sum_{j=1}^{r}a_{j}x_{j}]^{n}d\mu _{q}(x),
\end{equation*}%
where
\begin{equation*}
\int_{\mathbb{Z}_{p}^{r}}f(x)d\mu _{q}(x)=\int_{\mathbb{Z}_{p}}\int_{\mathbb{%
Z}_{p}}...\int_{\mathbb{Z}_{p}}f(x)d\mu _{q}(x_{1})d\mu _{q}(x_{2})...d\mu
_{q}(x_{r}).
\end{equation*}%
It is easily observed from (\ref{Eq-14}) that
\begin{eqnarray}
\beta _{n}^{(r)}(w &:&q\mid a_{1},...,a_{r})=\frac{1}{(1-q)^{n}}%
\sum_{l=0}^{n}\left(
\begin{array}{c}
n \\
l%
\end{array}%
\right) (-1)^{l}q^{wl}\prod_{j=1}^{r}\frac{(l+\frac{1}{a_{j}})}{[la_{j}+1]}
\label{Eq-15} \\
&=&\sum_{l=0}^{n}\left(
\begin{array}{c}
n \\
l%
\end{array}%
\right) [w]^{n-l}q^{wl}\beta _{l}^{(r)}(q\mid a_{1},...,a_{r}),  \notag
\end{eqnarray}%
for every positive integer $n$.

By (\ref{Eq-13}) and (\ref{Eq-15}), we note that
\begin{equation*}
\lim_{n\rightarrow \infty }\beta _{n}^{(r)}(w:q\mid
a_{1},...,a_{r})=B_{n}(w,r\mid a_{1},...,a_{r}).
\end{equation*}%
In the special case when%
\begin{equation*}
(a_{1},...,a_{r})=(1,...,1).
\end{equation*}%
we see that
\begin{equation*}
\beta _{n}^{(r)}(w:q\mid 1,1,...,1)=\beta _{n}^{(r)}(w:q),
\end{equation*}%
where $\beta _{n}^{(r)}(w:q)$ are the $q$-Bernoulli polynomials of order $r$
(see\cite{Kim5}), which reduces to the ordinary Bernoulli polynomials of
higher order $B_{n}^{(r)}(w)$ if $q=1$ (see for detail \cite{Kim8}, \cite%
{Kim10}, \cite{Kim11})

Kim\cite{Kim9} defined the analytic continuation of multiple zeta functions
( the Euler-Barnes multiple zeta functions ) depending on parameters $%
a_{1},...,a_{r}$ in the complex number field as follows:
\begin{equation}
\zeta _{r}(s,w,u\mid a_{1},...,a_{r})=\sum_{m_{1,...,m_{r}=0}}^{\infty
}u^{-(m_{1}+...+m_{r})}(w+m_{1}a_{1}+...+m_{N}a_{N})^{-s},  \label{eqq-20}
\end{equation}%
where $\Re (w)>0$ and $u\in \mathbb{C}$ with $\mid u\mid >1$ (\cite{Kim10}, %
\cite{Kim12}, \cite{Kim13}, \cite{Kim14}, \cite{Kim-Rim}).

We summarize our paper as follows:

In Section 2, of our paper, by using $q$-Volkenborn integration and uniform
differentiable function on $\mathbb{Z}_{p}$, we will construct $p$-adic $q$%
-zeta functions. This functions interpolate $q$-Bernoulli numbers. The
values of the $p$-adic $q$-zeta functions are given explicitly.

In Section 3, our primary aim is to give generating function of $q$%
-Bernoulli numbers and polynomials. These numbers and functions can be used
to prove analytic continuation of $q$-$L$-functions.

In Section 4, we give new generating functions which produce new definition
of Barnes' type Changhee $q$-Bernoulli polynomials and the generalized
Barnes' type Changhee $q$-Bernoulli numbers attached to Dirichlet character.
These functions are very important in constructing multiple zeta functions.
By using Mellin transformation formula, we also give relations between new
Barnes' type Changhee $q$-zeta functions and Barnes' type new Changhee $q$%
-Bernoulli numbers.

In Section 5, by using Mellin transformation formula of character generating
function of generalized Barnes' type Changhee $q$-Bernoulli numbers, we will
define the Dirichlet's type Changhee $q$-$L$-series. We give relations
between these functions, $q$-zeta functions and Dirichlet's type Changhee $q$%
-Bernoulli numbers, as well.

In Section 6, we will construct new generating function of multiple Changhee
$q$-Bernoulli polynomials. Under the Mellin transformation , we give
relation between this function and multiple $q$-zeta function, the multiple
Changhee $q$-Bernoulli numbers.

In Section 7, we give relation between $q$-gamma functions and zeta
functions. We also find some new results related to these functions.

In Section 8, we give analytic properties of $q$-$L$-function and $q$%
-Hurwitz zeta function. We \ prove relations between these functions and $q$%
-gamma function.

In Section 9, we define generalized multiple Changhee $q$-Bernoulli numbers
attached to the Drichlet character $\chi $. We also construct Dirichlet's
type multiple Changhee $q$-$L$-functions. These will lead to relations
between Dirichlet's type multiple Changhee $q$-$L$-functions and generalized
multiple Changhee $q$-Bernoulli numbers as well.

In Section 10, the main purpose is to prove analytic continuation of the
Euler-Barnes' type multiple $q$-Daehee zeta functions depending on the
parameters $a_{1},...,a_{r}$ which are taken positive parts in the Complex
Field. Thus, we construct generating function of $q$-Euler-Barnes' type
multiple Frobenius-Euler polynomials. We define Euler-Barnes' type multiple $%
q$-Daehee zeta Functions. Euler-Barnes' type multiple $q$-Daehee zeta
functions have a certain connection with Topology and Physics, together with
the algebraic relations among them. We give the values of these functions at
negative integers as well.

In Section 11, we give analytic continuation of Euler-Barnes' type Daehee $q$%
-zeta functions. We also give some remarks related to these functions.

\section{A Family of $p$-Adic $q$-Zeta Functions}

In this section, we need the following definitions and notations.

Every continuous function $f:$ $\mathbb{Z}_{p}\rightarrow \mathbb{K}$, $%
\mathbb{K}$ is an non-Archemedian valued field, has unique expansion as
\begin{equation}
f(x)=\sum_{n\geq 0}(\Delta _{q}^{n}f)(0)\left(
\begin{array}{c}
x \\
n%
\end{array}%
\right) _{q},  \label{EQ-1.1}
\end{equation}%
where
\begin{equation*}
(\Delta _{q}^{n}f)(0)\in \mathbb{K}
\end{equation*}
and
\begin{equation*}
(\Delta _{q}^{n}f)(0)\rightarrow 0\text{ as }n\rightarrow \infty
\end{equation*}
and $\Delta $ is the $q$-difference operator. We also have%
\begin{equation*}
\left(
\begin{array}{c}
m \\
n%
\end{array}%
\right) _{q}=\frac{[m]!}{[n]![m-n]!}\text{ \ ( }m\geq n\text{ )}
\end{equation*}
( see, for example \cite{T. Kim-SD.Kim-DW.Park}).

\begin{definition}
A function $f$ is called uniform differentiable function if it satisfies the
following conditions:

1) $F_{f}:\mathbb{Z}_{p}\times \mathbb{Z}_{p}\rightarrow \mathbb{C}_{p}$,

$F_{f}(x,y)=\frac{f(x)-f(y)}{x-y}.$

2) $\lim_{x\rightarrow y}F_{f}(x,y)=f^{^{\prime }}(y).$
\end{definition}

From (\ref{EQ-1.1}), let $f:$ $\mathbb{Z}_{p}\rightarrow \mathbb{Z}_{p}$ be
a function. Then $f$ can be written as
\begin{equation}
f(x)=\sum_{n=0}^{\infty }\left(
\begin{array}{c}
x \\
n%
\end{array}
\right) _{q}\Delta ^{n}f(0),  \label{EQ-1.2}
\end{equation}
where
\begin{equation*}
\Delta ^{n}f(0)=\sum_{i=0}^{n}\left(
\begin{array}{c}
n \\
i%
\end{array}
\right) _{q}(-1)^{n-i}f(i).
\end{equation*}
By using (\ref{EQ-1.2}), we have

1) $f$ is continuous function $\Leftrightarrow $ $\mid \Delta ^{n}f(0)\mid
_{p}\rightarrow 0$ as $n\rightarrow \infty $ .

2) $f$ is differentiable function $\Leftrightarrow $ $\mid \frac{\Delta
^{n}f(0)}{n}\mid _{p}\rightarrow 0$ as $n\rightarrow \infty $ .

3) $f$ is analytic function $\Leftrightarrow $ $\mid \frac{\Delta ^{n}f(0)}{%
n!}\mid _{p}\rightarrow 0$ as $n\rightarrow \infty $ .

4) $f$ is uniform differentiable function, $f\in UD$ $\Leftrightarrow n\mid
\Delta ^{n}f(0)\mid _{p}\rightarrow 0$ as $n\rightarrow \infty $ .

5) $\ f\in C^{(n)}$ function $\Leftrightarrow \lim_{m\rightarrow \infty
}m^{n}\mid \Delta ^{m}f(0)\mid _{p}\rightarrow 0$ ,

where%
\begin{equation*}
C^{(m)}=C^{(m)}(\mathbb{Z}_{p},\mathbb{K)=}\left\{ f:f:\mathbb{Z}%
_{p}\rightarrow \mathbb{K}\text{, }m\text{ times strictly differentiable}%
\right\} .
\end{equation*}%
If $f\in C^{(m)}(\mathbb{Z}_{p},\mathbb{K)}$, then $f^{(m-1)}\in C^{(1)}$.
We also note that $(\Delta _{q}^{n}f)(0)$ is the $n$th Mahler coefficient of
$f$ at $q=1$.

The function $f$ is differentiable at%
\begin{equation*}
x\in \mathbb{Z}_{p}\Leftrightarrow \lim_{m\rightarrow \infty }\frac{(\Delta
_{q}^{m}f)(0)}{[m]}=0,
\end{equation*}

in \ which case
\begin{equation*}
f^{(1)}(x)=\frac{\log q}{q-1}\sum_{m=1}^{\infty }\frac{(\Delta _{q}^{m}f)(x)%
}{[m]}(-1)^{m-1}q^{-\left(
\begin{array}{c}
m \\
2%
\end{array}%
\right) },
\end{equation*}%
where $f^{(1)}$\ denotes first derivative. We say that \ $f$ is strictly
differentiable at a point $a\in \mathbb{Z}_{p}$, and denote this property by
$f\in C^{(1)}$, if the difference quotients
\begin{equation*}
F_{f}(x,y)=\frac{f(x)-f(y)}{x-y},
\end{equation*}%
have a limit $l=f^{^{\prime }}(a)$ as $(x,y)\rightarrow (a,a)$. Recall that
if $f\in C^{(1)}(\mathbb{Z}_{p},\mathbb{C}_{p}\mathbb{)}$ function, then
there exist a unique continuos function
\begin{equation*}
(sf)(x)=\sum_{k\func{mod}x}f(k).
\end{equation*}
The function $sf$ satisfies the following properties:

1) $\ (sf)(x+1)-(sf)(x)=f(x)$,

2) If $f\in C^{(1)}(\mathbb{Z}_{p},\mathbb{C}_{p}\mathbb{)}$, then
\begin{equation*}
(sf)\in C^{(1)}(\mathbb{Z}_{p},\mathbb{C}_{p}\mathbb{)},
\end{equation*}%
and
\begin{equation*}
\parallel sf\parallel _{1}\leq p\parallel f\parallel _{1},
\end{equation*}%
where
\begin{equation*}
\parallel f\parallel _{1}=\parallel f\parallel _{p}\vee \parallel \Delta
_{n}f\parallel _{\infty }
\end{equation*}%
and (cf. \cite{T. Kim-SD.Kim-DW.Park})%
\begin{equation*}
\parallel f\parallel _{m}=\mid f^{(m-1)}\mid _{\sup }\vee \mid
F_{f^{(m-1)}}(x,y)\mid _{\sup },
\end{equation*}
\begin{equation*}
\parallel f\parallel _{\infty }=\sup_{x\in \mathbb{Z}_{p}}\mid f(x)\mid _{p}.
\end{equation*}%
We note that
\begin{equation*}
\Delta _{n}f(m_{1},...,m_{n};x)=\Delta _{1}(\Delta
_{n-1}f(m_{1},...,m_{n-1};x)),
\end{equation*}%
and

\begin{equation*}
\Delta _{1}f(m,x)=\frac{f(x+m)-f(x)}{n}.
\end{equation*}

Therefore, by (\ref{Eq-17}), it is easy to see that
\begin{equation*}
\int_{\mathbb{Z}_{p}}q^{-x}f(x)d\mu _{q}(x)=\lim_{N\rightarrow \infty }\frac{%
sf(p^{N})-sf(0)}{[p^{N}]}=\frac{q-1}{\log q}(sf)^{^{\prime }}(0),
\end{equation*}%
(see \cite{T. Kim-SD.Kim-DW.Park}).

\begin{proposition}
1) If $f\in C^{(1)}(\mathbb{Z}_{p},\mathbb{C}_{p}\mathbb{)}$, then
\begin{equation*}
\int_{\mathbb{Z}_{p}}q^{-x}f(x)d\mu _{q}(x)\leq p\parallel f\parallel _{1}.
\end{equation*}%
that is, the $q$-Volkenborn integral a linear, continuos function on $%
C^{(1)}(\mathbb{Z}_{p},\mathbb{C}_{p}\mathbb{).}$

2) Let
\begin{equation*}
f(x)=\sum_{n=0}^{\infty }\left(
\begin{array}{c}
x \\
n%
\end{array}%
\right) _{q}\Delta ^{n}f(0).
\end{equation*}%
Then the $q$-Mahler representation of the $C^{(1)}$-function $f$ satisfies
the following equation:
\begin{equation*}
\int_{\mathbb{Z}_{p}}f(x)d\mu _{q}(x)=\sum_{n=0}^{\infty }(\Delta ^{n}f)(0)%
\frac{(-1)^{n}}{[n+1]}q^{n+1-\left(
\begin{array}{c}
n+1 \\
2%
\end{array}%
\right) }.
\end{equation*}%
3) In particular, the Fourier transformation of the $q$-Volkenborn integral
is given by
\begin{equation}
\int_{\mathbb{Z}_{p}}q^{-x}T^{[x]}d\mu _{q}(x)=\frac{q-1}{\log q}T^{\frac{1}{%
1-q}}-\log T\sum_{n=0}^{\infty }q^{n}T^{[n]}.  \label{EQ-1.3}
\end{equation}

4)
\begin{equation}
\int_{\mathbb{Z}_{p}}q^{-x}[x]^{n}d\mu _{q}(x)=\beta _{n}(q),  \label{EQ-1.4}
\end{equation}%
where $\beta _{n}(q)$ the $n$th $q$-Bernoulli numbers (see, for detail \cite%
{Kim6}).

5) If%
\begin{equation*}
f\in C^{(1)}(\mathbb{Z}_{p},\mathbb{C}_{p}),
\end{equation*}%
$n\geq 0$, $j\in \left\{ 0,1,2,...,p-1\right\} $ and%
\begin{equation*}
\mathbb{T}_{p}=\mathbb{Z}_{p}\setminus p\mathbb{Z}_{p},
\end{equation*}%
then
\begin{eqnarray*}
\int_{j+p^{n}\mathbb{Z}_{p}}q^{-(j+xp^{n})}f(x)d\mu _{q}(x)
&=&\lim_{N\rightarrow \infty }\sum_{x=0}^{p^{N}-1}f(j+p^{n}x)\frac{1}{%
[p^{n+N}]} \\
&=&\frac{1}{[p^{n}]}\lim_{N\rightarrow \infty }\frac{1}{[p^{N}:q^{p^{n}}]}%
\sum_{x=0}^{p^{N}-1}f(j+p^{n}x) \\
&=&\frac{1}{[p^{n}]}\int_{\mathbb{Z}_{p}}q^{-xp^{n}}f(j+p^{n}x)d\mu
_{q^{p^{n}}}(x),
\end{eqnarray*}%
and
\begin{eqnarray*}
\int_{\mathbb{T}_{p}}f(x)d\mu _{q}(x) &=&\int_{\mathbb{Z}_{p}\setminus p%
\mathbb{Z}_{p}}f(x)d\mu _{q}(x) \\
&=&\int_{\mathbb{Z}_{p}}f(x)d\mu _{q}(x)-\frac{1}{[p]}\int_{\mathbb{Z}%
_{p}}f(px)d\mu _{q^{p}}(x) \\
&=&\frac{1}{[p]}\left( [p]\int_{\mathbb{Z}_{p}}f(x)d\mu _{q}(x)-\int_{%
\mathbb{Z}_{p}}f(px)d\mu _{q^{p}}(x)\right) .
\end{eqnarray*}
\end{proposition}

We note that, substituting $T=e^{t}$ into (\ref{EQ-1.3}), than we obtain (%
\ref{EQ-1.4}).

We now consider $p$-adic $q$-zeta function
\begin{equation*}
\zeta _{p,q,1},\zeta _{p,q,2},...,\zeta _{p,q,p-1},
\end{equation*}%
which are given by
\begin{equation*}
\zeta _{p,q,j}(s)=\frac{1}{j+(p-1)s}\int_{\mathbb{T}%
_{p}}q^{-x}[x]^{j}[x]^{(p-1)s}d\mu _{q}(x),\text{ }(j=0,1,2,...,p-1),
\end{equation*}%
where%
\begin{equation*}
\mathbb{T}_{p}=\mathbb{Z}_{p}\setminus p\mathbb{Z}_{p}\text{, and}\mid s\mid
_{p}<p^{\frac{p-2}{p-1}}\text{ \ ( }s\neq -\frac{j}{p-1}\text{).}
\end{equation*}

We note that
\begin{eqnarray*}
\int_{\mathbb{Z}_{p}\setminus p\mathbb{Z}_{p}}q^{-x}[x]^{j}[x]^{(p-1)s}d\mu
_{q}(x) &=&\int_{\mathbb{T}_{p}}q^{-x}[x]^{j}[x]^{(p-1)s}d\mu _{q}(x) \\
&=&\int_{\mathbb{T}_{p}}q^{-x}[x]^{j}\sum_{n=0}^{\infty }s^{n}(p-1)^{n}\frac{%
(\log [x])^{n}}{n!}d\mu _{q}(x) \\
&=&\sum_{n=0}^{\infty }a_{n,q}s^{n},
\end{eqnarray*}%
where, as before,%
\begin{equation*}
\mid s\mid _{p}<p^{\frac{p-2}{p-1}}\text{ \ ( }s\neq -\frac{j}{p-1}\text{),}
\end{equation*}%
and
\begin{equation*}
a_{n,q}=\frac{(p-1)^{n}}{n!}\int_{\mathbb{T}_{p}}q^{-x}[x]^{j}\log
^{n}[x]d\mu _{q}(x).
\end{equation*}%
If $n\rightarrow \infty $, then $\mid a_{n,q}\mid _{p}\rightarrow 0$. By
using the above definition, we obtain
\begin{eqnarray*}
\zeta _{p,q,0}(s) &=&\frac{1}{(p-1)s}\int_{\mathbb{T}%
_{p}}q^{-x}[x]^{(p-1)s}d\mu _{q}(x) \\
&=&\frac{1}{(p-1)s}\sum_{n=0}^{\infty }a_{n,q}s^{n},
\end{eqnarray*}%
where
\begin{eqnarray*}
a_{0,q} &=&\int_{\mathbb{T}_{p}}q^{-x}d\mu _{q}(x)=\int_{\mathbb{Z}%
_{p}}q^{-x}d\mu _{q}(x)-\int_{p\mathbb{Z}_{p}}q^{-x}d\mu _{q}(x) \\
&=&\frac{q-1}{\log q}-\frac{q-1}{\log q}\frac{1}{p}=\frac{(q-1)(p-1)}{p\log q%
}.
\end{eqnarray*}

Therefore, we arrive at the following result.

\begin{theorem}
By means of the\textit{\ following transformation: }%
\begin{equation*}
s\rightarrow \zeta _{p,q,0}-\frac{1}{ps}\left( \frac{q-1}{\log q}\right) ,
\end{equation*}%
each of the functions
\begin{equation*}
\zeta _{p,q,1},\text{ }\zeta _{p,q,2},...,\zeta _{p,q,p-1}
\end{equation*}%
can be extend to the corresponding analytic function on the following set:
\begin{equation*}
\mathit{\ }\mathbb{B}\mathit{:=}\left\{ s:s\in \mathbb{C}_{p}\text{ and}\mid
s\mid _{p}<p^{\frac{p-2}{p-1}}\right\} .
\end{equation*}
\end{theorem}

\begin{remark}
It easily follows from the above observations that%
\begin{equation*}
\zeta _{p,q,0}(s)=\frac{q-1}{\log q}\frac{1}{ps}+\Theta (s),
\end{equation*}%
where $\Theta (s)$ is analytic function. Thus, before giving the connection
between the $p$-adic $q$-zeta functions and the classical $q$-zeta
functions, we determine the values of%
\begin{equation*}
\zeta _{p,q,1}(s),\text{ }\zeta _{p,q,2}(s),...,\zeta _{p,q,p-1}(s).
\end{equation*}
\ Then, by using these values as well as the $p$-adic interpolation of
sequences of values of $\zeta _{q}$ at certain negative integers , we will
construct the $p$-adic $q$-zeta functions.
\end{remark}

The following theorems provides us with the relationship between $\zeta
_{p,q,j}(s)$ and $\beta _{n}(q)$ are given by

\begin{proposition}
Let $n$ and $p$ be positive integers with $p$ prime. Then
\begin{equation*}
\zeta _{p,q,0}(n)=\frac{1}{(p-1)n}(\beta _{(p-1)n}(q)-[p]^{(p-1)n-1}\beta
_{(p-1)n-1}(q^{p})).
\end{equation*}
\end{proposition}

\begin{proof}
\begin{eqnarray*}
\zeta _{p,q,0}(s) &=&\frac{1}{(p-1)s}\int_{T_{p}}q^{-x}[x]^{(p-1)s}d\mu
_{q}(x) \\
&=&\frac{1}{(p-1)s}\int_{\mathbb{Z}_{p}\setminus p\mathbb{Z}%
_{p}}q^{-x}[x]^{(p-1)s}d\mu _{q}(x) \\
&=&\frac{1}{(p-1)s}\left( \int_{\mathbb{Z}_{p}}q^{-x}[x]^{(p-1)s}d\mu
_{q}(x)-\int_{p\mathbb{Z}_{p}}q^{-x}[x]^{(p-1)s}d\mu _{q}(x)\right) \\
&=&\frac{1}{(p-1)s}\left( \int_{\mathbb{Z}_{p}}q^{-x}[x]^{(p-1)s}d_{q}x-%
\int_{\mathbb{Z}_{p}}q^{-px}[px]^{(p-1)s}d\mu _{q^{p}}(x)\right)
\end{eqnarray*}%
By setting $s=n$, ( $n\in \mathbb{Z}^{+}$ ) and using (\ref{EQ-1.4}) in the
above, we easily obtain
\begin{equation*}
\zeta _{p,q,0}(n)=\frac{1}{(p-1)n}\left( \beta
_{(p-1)n}(q)-[p]^{(p-1)n-1}\beta _{(p-1)n-1}(q^{p})\right) ,
\end{equation*}%
which completes our proof of Proposition 2.
\end{proof}

By applying a similar method for $\zeta _{p,q,j}$, we arrive at the desired
result asserted by Theorem 2 below.

\begin{theorem}
Let $n,p,j$ be positive integers with $p$ prime. Then
\begin{equation*}
\zeta _{p,q,j}(n)=\frac{1}{j+(p-1)n}(\beta
_{j+(p-1)n}(q)-[p]^{j+(p-1)n-1}\beta _{j+(p-1)n-1}(q^{p})).
\end{equation*}
\end{theorem}

The proof of Theorem 2 is simillar to that of Proposition 2. So we omit it.

The classical $q$-zeta function was defined by Kim\cite{Kim-Rim} as follows:
\begin{equation*}
\zeta _{q}(s)=\sum_{n=0}^{\infty }\frac{q^{n}}{[n]^{s}}-\frac{1}{s-1}\frac{%
(1-q)^{s}}{\log q},
\end{equation*}
for $s\in \mathbb{C.}$

For any positive integer $n$,\ we have
\begin{equation*}
\zeta _{q}(1-n)=-\frac{\beta _{n}(q)}{n}.
\end{equation*}%
Furthermore, if we define
\begin{equation*}
\zeta _{q}^{\ast }(s)=\zeta _{q}(s)-[p]^{-s}\zeta _{q^{p}}(s),
\end{equation*}%
then, for $j\in \left\{ 0,1,2,...,p-1\right\} $, and $n\geq 0$, we find that
\begin{equation*}
\zeta _{p,q,j}(n)=-\zeta _{q}^{\ast }(1-(j+(p-1)n).
\end{equation*}

\section{The $q$-Bernoulli Numbers and the $q$-Bernoulli Polynomials}

Our primary aim in this section is to give generating functions of $q$%
-Bernoulli numbers and $q$-Bernoulli polynomials. These numbers will be used
to prove analytic continuation of $q$-$L$-series.

We first define $q$-version of each of the functions $F(t,x)$ and $F(t)$
occurring in (\ref{eq-1}) and (\ref{q-2}), respectively. The generating
function $F_{q}(t)$ of $q$-Bernoulli numbers $\beta _{n}(q)$ ( $n\geq 0$ )
is given by \cite{Kim7}:
\begin{equation}
F_{q}(t)=\frac{q-1}{\log q}\exp (\frac{t}{1-q})-t\sum_{n=0}^{\infty
}q^{n}e^{[n]t}=\sum_{n=0}^{\infty }\frac{\beta _{n}(q)t^{n}}{n!},
\label{EQ-2.1}
\end{equation}%
We note that the definition (\ref{q-2}) and (\ref{EQ-2.1}) that

\begin{equation*}
\lim_{q\rightarrow 1}\beta _{n}(q)=B_{n},
\end{equation*}%
and
\begin{equation*}
\lim_{q\rightarrow 1}F_{q}(t)=F(t)=\frac{t}{e^{t}-1}.
\end{equation*}%
The generating function $F_{q}(x,t)$ of the $q$-Bernoulli polynomials $\beta
_{n}(x:q)$ ( $n\geq 0$ ) is defined analogously as follows:
\begin{equation}
F_{q}(x,t)=\frac{q-1}{\log q}\exp (\frac{t}{1-q})-t\sum_{n=0}^{\infty
}q^{n+x}e^{[n+x]t}=\sum_{n=0}^{\infty }\frac{\beta _{n}(x:q)t^{n}}{n!}.
\label{EQ-2.2}
\end{equation}%
We note from (\ref{eq-1}) and (\ref{EQ-2.2}) that

\begin{equation*}
\lim_{q\rightarrow 1}\beta _{n}(x:q)=B_{n}(x),
\end{equation*}%
and
\begin{equation*}
\lim_{q\rightarrow 1}F_{q}(x,t)=F(x,t)=\frac{te^{xt}}{e^{t}-1}.
\end{equation*}

The remarkable point here is that the series on the righet-hand side of (\ref%
{EQ-2.1}) and (\ref{EQ-2.2}) are uniformly convergent in the wider sense.
Therefore, we shall explicitly determine the $q$-Bernoulli numbers as
follows:
\begin{equation*}
\beta _{0,q}=\frac{q-1}{\log q},\text{ \ }q(q\beta _{q}+1)^{n}-\beta
_{n}(q)=\left\{
\begin{array}{c}
1,\text{ if }n=1 \\
0,\text{ if }n>1,%
\end{array}%
\right.
\end{equation*}%
with the usual convention about replacing $\beta ^{n}$ by $\beta _{n}$.

For the $q$-Bernoulli polynomials are defined by\ (\ref{EQ-2.2}), we first
derive an explicit representation given by Theorem 3 below.

\begin{theorem}
\begin{equation*}
\beta _{n}(x:q)=\sum_{l=0}^{n}\left(
\begin{array}{c}
n \\
l%
\end{array}%
\right) q^{lx}\beta _{l}(q)[x]^{n-l}.
\end{equation*}
\end{theorem}

\begin{proof}
By using Cauchy product in (\ref{EQ-2.1}) and (\ref{EQ-2.2}), we have
\begin{equation*}
\sum_{n=0}^{\infty }\frac{\beta _{n}(x:q)t^{n}}{n!}=\sum_{n=0}^{\infty
}(\sum_{l=0}^{n}\left(
\begin{array}{c}
n \\
l%
\end{array}%
\right) q^{lx}\beta _{l}(q)[x]^{n-l})\frac{t^{n}}{n!}.
\end{equation*}%
After some elementary calculations in the above, we easily arrive at the
desired result.
\end{proof}

We note that
\begin{eqnarray*}
\beta _{n}(x &:&q)=\sum_{l=0}^{n}\left(
\begin{array}{c}
n \\
l%
\end{array}%
\right) q^{lx}\beta _{l}(q)[x]^{n-l} \\
&=&(q^{x}\beta (q)+[x])^{n}.
\end{eqnarray*}%
Now we construct generalized $q$-Bernoulli numbers associated with a
Dirichlet characater.

Let $\chi $ be a Dirichlet character of conductor $f\in \mathbb{Z}^{+}$. \
We define the generating function of generalized $q$-Bernoulli numbers
attached to $\chi $ as follows:
\begin{eqnarray}
F_{q,\chi }(t) &=&-t\sum_{a=1}^{f}\chi (a)\sum_{n=0}^{\infty
}q^{fn+a}e^{[fn+a]t}  \label{EQ-2.3} \\
&=&-t\sum_{n=0}^{\infty }\chi (n)q^{n}e^{[n]t}  \notag \\
&=&\sum_{n=0}^{\infty }\beta _{n,\chi }(q)\frac{t^{n}}{n!}.  \notag
\end{eqnarray}%
where the coefficients, $\beta _{n,\chi }(q)$ ( $n\geq 0$ ) are called
generalized $q$-Bernoulli numberswith a Dirichlet characater. We note from
the definitions in (\ref{eq-3}) and (\ref{EQ-2.3}) that

\begin{equation*}
\lim_{q\rightarrow 1}\beta _{n,\chi }(q)=B_{n,\chi },
\end{equation*}%
and
\begin{equation*}
\lim_{q\rightarrow 1}F_{q,\chi }(t)=F_{\chi }(t)=\sum_{a=1}^{f}\chi (a)\frac{%
te^{at}}{e^{tf}-1}=\sum_{n=0}^{\infty }B_{n,\chi }\frac{t^{n}}{n!}.
\end{equation*}

By using (\ref{EQ-2.3}), we also have
\begin{equation}
F_{q,\chi }(t)=\frac{1}{[f]}\sum_{a=1}^{f}\chi (a)F_{q^{f}}(\frac{a}{f},[f]t)
\label{EQ-2.5}
\end{equation}

By applying Cauchy product in (\ref{EQ-2.2}), (\ref{EQ-2.3}) and (\ref%
{EQ-2.5}), we easily obtain the following theorems:

\begin{theorem}
Let $\chi $ be a Dirichlet character with conductor $f\in \mathbb{Z}^{+}$.
Then we have
\begin{equation*}
\beta _{n,\chi }(q)=[f]^{n-1}\sum_{a=1}^{f}\chi (a)\beta _{n}(\frac{a}{f}%
:q^{f}).
\end{equation*}
\end{theorem}

We now construct generating function of generalized $q$-Bernoulli
polynomials associated with a Dirichlet characater as follows:
\begin{eqnarray}
F_{q,\chi }(x,t) &=&q^{x}te^{-[x]t}\sum_{n=0}^{\infty }\chi
(n)q^{n}e^{[n]q^{x}t}  \label{EQ-2.6} \\
&=&-t\sum_{n=0}^{\infty }\chi (n)q^{n+x}e^{[n+x]t}  \notag \\
&=&\sum_{n=0}^{\infty }\beta _{n,\chi }(x:q)\frac{t^{n}}{n!}.  \notag
\end{eqnarray}%
By using (\ref{EQ-2.2}), (\ref{EQ-2.3}) and (\ref{EQ-2.6}), we obtain
\begin{eqnarray*}
F_{q,\chi }(x,t) &=&-t\sum_{n=0}^{\infty }\chi (n)q^{n+x}e^{[n+x]t} \\
&=&-e^{[x]t}q^{x}t\sum_{a=1}^{f}\chi (a)\sum_{n=0}^{\infty
}q^{fn+a}e^{[fn+a]tq^{x}},
\end{eqnarray*}%
which, in view of the following well-known identity
\begin{equation*}
\lbrack x+a]=[x]+q^{x}[a]
\end{equation*}%
yields
\begin{equation*}
F_{q,\chi }(x,t)=-t\sum_{a=1}^{f}\chi (a)\sum_{n=0}^{\infty
}q^{fn+a+x}e^{[fn+a+x]t}.
\end{equation*}%
After some elementary calculation we arrive at the following theorem:

\begin{theorem}
Let $\chi $ be a Dirichlet character of conductor $f\in \mathbb{Z}^{+}$.
Then we have
\begin{equation}
F_{q,\chi }(x,t)=\frac{1}{[f]}\sum_{a=1}^{f}\chi (a)F_{q^{f}}(\frac{a+x}{f}%
,[f]t)  \label{EQ-2.7}
\end{equation}
\end{theorem}

Note that substituting $x=0$ into (\ref{EQ-2.7}), then we obtain (\ref%
{EQ-2.5}).

By comparing the coefficients on both sides of (\ref{EQ-2.6}) and (\ref%
{EQ-2.7}), we easily see that
\begin{eqnarray*}
\beta _{n,\chi }(x &:&q)=\sum_{l=0}^{n}\left(
\begin{array}{c}
n \\
l%
\end{array}%
\right) q^{lx}[x]^{n-l}\sum_{a=1}^{f}\chi (a)\beta _{l}(\frac{a}{f}:q^{f}) \\
&=&\sum_{l=0}^{n}\left(
\begin{array}{c}
n \\
l%
\end{array}%
\right) q^{lx}[x]^{n-l}\beta _{l,\chi }(q).
\end{eqnarray*}%
By using definition of $\beta _{l,\chi }(q)$ into the above, we have obtain
the following result.

\begin{theorem}
Let $\chi $ be a Dirichlet character of conductor $f\in \mathbb{Z}^{+}$.
Then
\begin{equation*}
\beta _{n,\chi }(x:q)=\frac{1}{[f]^{1-n}}\sum_{a=0}^{f-1}\chi (a)\beta _{n}(%
\frac{a+x}{f}:q^{f}).
\end{equation*}
\end{theorem}

\section{A Class of $q$-Multiple Zeta Functions}

In this section, we give new generating functions which produce new
definitions of Barnes' type of Changhee $q$-Bernoulli polynomials and the
generalized Barnes' type Changhee $q$-Bernoulli numbers with attached to $%
\chi $, Dirichlet character with conductor with conductor $f\in \mathbb{Z}%
^{+}$. These generating functions are very important in case of multiple
zeta function. Therefore, by using these generating functions, we will give
relation between Barnes' type Changhee $q$-zeta function and Barnes' type
Changhee $q$-Bernoulli numbers.

Let $w,w_{1},w_{2},\ldots ,w_{r}$ be complex numbers such that $w_{i}\neq 0$
for $i=1,2,\ldots ,r$. We define Barnes' type of Changhee $q$-Bernoulli
polynomials of $w$ with parameters $w_{1}$ as follows:
\begin{eqnarray}
F_{q}(w,t &\mid &w_{1})=\frac{q-1}{\log q}e^{\frac{t}{1-q}%
}-w_{1}t\sum_{n=0}^{\infty }q^{w_{1}n+w}e^{[w_{1}n+w]t}  \label{EQ-3.1} \\
&=&\sum_{n=0}^{\infty }\frac{\beta _{n}(w:q\mid w_{1})t^{n}}{n!}\text{ \ (}%
\mid t\mid <2\pi \text{ ),}  \notag
\end{eqnarray}%
where the coefficients, $\beta _{n}(w:q\mid w_{1})$ ( $n\geq 0$ ) are called
Barnes' type of Changhee $q$-Bernoulli polynomials in $w$ with parameters $%
w_{1}$.

We note that

\begin{equation*}
\lim_{q\rightarrow 1}\beta _{n}(w:q\mid w_{1})=w_{1}^{n}\beta _{n}(w),
\end{equation*}%
and
\begin{equation*}
\lim_{q\rightarrow 1}F_{q}(w,t\mid w_{1})=\frac{w_{1}t}{e^{w_{1}t}-1}e^{wt},
\end{equation*}%
where $\beta _{n}(w)$ are the ordinary Barnes Bernoulli polynomials.

By using (\ref{EQ-3.1}), we easily obtain\cite{Kim-Rim}, \cite{Kim10}
\begin{equation*}
\beta _{n}(w:q\mid w_{1})=\frac{1}{(1-q)^{n}}\sum_{l=0}^{n}\left(
\begin{array}{c}
n \\
l%
\end{array}%
\right) q^{lw}(-1)^{l}\frac{lw_{1}}{[lw_{1}]}.
\end{equation*}%
If $w=0$ in the above, then%
\begin{equation*}
\beta _{n}(0:q\mid w_{1})=\beta _{n}(q\mid w_{1}),
\end{equation*}%
where $\beta _{n}(q\mid w_{1})$ are called Barnes' type Changhee $q$%
-Bernoulli numbers with parameter $w_{1}$.

Let $\chi $ be the Dirichlet character with conductor $f$. Then the
generalized Barnes' type Changhee $q$-Bernoulli numbers with attached to $%
\chi $ are defined as follows:
\begin{eqnarray}
F_{q,\chi }(t &\mid &w_{1})=-w_{1}t\sum_{n=1}^{\infty }\chi
(n)q^{w_{1}n}e^{[w_{1}n]t}  \label{EQ-3.2} \\
&=&\sum_{n=0}^{\infty }\frac{\beta _{n,\chi }(q\mid w_{1})t^{n}}{n!}\text{ \
(}\mid t\mid <2\pi \text{ ).}  \notag
\end{eqnarray}%
We easily see from (\ref{EQ-3.2}) that
\begin{equation}
F_{q,\chi }(t\mid w_{1})=\frac{1}{[f]}\sum_{a=1}^{f}\chi (a)F_{q^{f}}(\frac{%
w_{1}a}{f},[f]t\mid w_{1}).  \label{EQ-3.3}
\end{equation}%
Now by using (\ref{EQ-3.1}), (\ref{EQ-3.2}) and (\ref{EQ-3.3}), and after
some elementary calculations, we arrive at the following theorem.

\begin{theorem}
Let $\chi $ be a Dirichlet character of conductor $f\in \mathbb{Z}^{+}$.
Then
\begin{equation*}
\beta _{n,\chi }(q\mid w_{1})=\frac{1}{[f]^{1-n}}\sum_{a=0}^{f-1}\chi
(a)\beta _{n}(\frac{aw_{1}}{f}:q^{f}\mid w_{1}).
\end{equation*}
\end{theorem}

By applying Mellin transformation in (\ref{EQ-3.1}), we obtain
\begin{equation}
\frac{1}{\Gamma (s)}\int_{0}^{\infty }t^{s-2}F_{q}(w,-t\mid w_{1})dt=-\frac{%
(1-q)^{s}}{s-1}\frac{1}{\log q}+w_{1}\sum_{n=0}^{\infty }\frac{q^{w_{1}n+w}}{%
[w_{1}n+w]^{s}},  \label{EQ-3.4}
\end{equation}
where $\Gamma (s)$ is denoted Euler gamma function.

Note that, by substituting $w=w_{1}=q=1$ into (\ref{EQ-3.4}), then we obtain
Hurwitz zeta function.

We define Barnes' type Changhee $q$-zeta function as follows:

\begin{definition}
For $s\in \mathbb{C}$, we have
\begin{equation}
\zeta _{q}(s,w\mid w_{1})=-\frac{(1-q)^{s}}{s-1}\frac{1}{\log q}%
+w_{1}\sum_{n=0}^{\infty }\frac{q^{w_{1}n+w}}{[w_{1}n+w]^{s}}.
\label{EQ-3.5}
\end{equation}
\end{definition}

\begin{theorem}
If $n\in \mathbb{Z}^{+}$, then
\begin{equation*}
\zeta _{q}(1-n,w\mid w_{1})=-\frac{\beta _{n}(w:q\mid w_{1})}{n}.
\end{equation*}
\end{theorem}

\begin{proof}
In view of (\ref{EQ-3.4}), we define $Y(s)$ by means of the following
contour integral:%
\begin{equation}
Y(s)=\int_{C}z^{s-2}F_{q}(w,-z\mid w_{1})dz,  \label{EQ-355}
\end{equation}%
where $C$ is Hankel's contour along the cut joining the points $z=0$ and $%
z=\infty $ on the real axis, which starts from the point at $\infty $,
encircles the origin ( $z=0$ ) once in the positive (counter-clockwise)
direction, and returns to the point at $\infty $ ( see for details, \cite{E.
T. Wittaker and G. N. Watson} p. 245). Here, as usual,we interpret $z^{s}$
to mean $\exp (s\log z)$, where we asume $\log $ to defined by $\log t$ on
the top part of the real axis and by $\log t+2\pi i$ on the bottom part of
the real axis. We thus find from the definition (\ref{EQ-355}) that
\begin{equation*}
Y(s)=(e^{2\pi is}-1)\int_{\varepsilon }^{\infty }t^{s-2}F_{q}(w,-t\mid
w_{1})dt
\end{equation*}%
\begin{equation*}
+\int_{C_{\varepsilon }}z^{s-2}F_{q}(w,-z\mid w_{1})dz,
\end{equation*}%
where $C_{\varepsilon }$ denotes a circle of radius $\varepsilon >0$ (and
centred at the origin), which is described in the positive
(counter-clockwise) direction. Assume first that $\func{Re}(s)>1$. Then
\begin{equation*}
\int_{C_{\varepsilon }}\rightarrow 0\text{ as }\varepsilon \rightarrow 0,
\end{equation*}%
so we have
\begin{equation*}
Y(s)=(e^{2\pi is}-1)\int_{0}^{\infty }t^{s-2}F_{q}(w,-t\mid w_{1})dt,
\end{equation*}%
which, upon substituting from (\ref{EQ-3.1}) into it, yields%
\begin{eqnarray*}
Y(s) &=&(e^{2\pi is}-1)(-\frac{1-q}{\log q}\int_{0}^{\infty }t^{s-2}\exp (-%
\frac{t}{1-q})dt \\
&&+w_{1}\sum_{n=0}^{\infty }q^{w_{1}n+w}\int_{0}^{\infty
}t^{s-1}e^{-[w_{1}n+w]t}dt).
\end{eqnarray*}%
After some elementary calculations, we thus find that%
\begin{equation*}
Y(s)=(e^{2\pi is}-1)\Gamma (s)\left( -\frac{(1-q)^{s}}{(s-1)\log q}%
+w_{1}\sum_{n=0}^{\infty }\frac{q^{w_{1}n+w}}{[w_{1}n+w]^{s}}\right) .
\end{equation*}%
By using (\ref{EQ-3.5}) in the above, we get%
\begin{equation*}
Y(s)=(e^{2\pi is}-1)\Gamma (s)\zeta _{q}(s,w\mid w_{1}).
\end{equation*}%
Therefore%
\begin{equation}
\zeta _{q}(s,w\mid w_{1})=\frac{Y(s)}{(e^{2\pi is}-1)\Gamma (s)},
\label{EQ-356}
\end{equation}%
which, by analytic continuation, holds true for all $s\neq 1$. This
evidently provides us with an analytic continuation of $\zeta _{q}(s,w\mid
w_{1})$.

We now consider the situation when we let $s\rightarrow 1-n$ in (\ref{EQ-356}%
), where $n$ is a positive integer. Then since%
\begin{equation*}
e^{2\pi is}=e^{2\pi i(1-n)}=1\text{ \ ( }n\in \mathbb{Z}^{+}\text{ ),}
\end{equation*}%
we have the following limit relationship:%
\begin{eqnarray}
\lim_{s\rightarrow 1-n}\left\{ (e^{2\pi is}-1)\Gamma (s)\right\}
&=&\lim_{s\rightarrow 1-n}\left\{ \frac{(e^{2\pi is}-1)}{\sin (\pi s)}\frac{%
\pi }{\Gamma (1-s)}\right\}  \label{EQ-358} \\
&=&\frac{2\pi i(-1)^{n-1}}{(n-1)!}\text{ \ ( }n\in \mathbb{Z}^{+}\text{ )}
\notag
\end{eqnarray}%
by means of the familiar reflection formula for $\Gamma (s)$. Furthermore,
since the integrand in (\ref{EQ-355}) has simple pole order $n+1$ at $z=0$,
where also find from the definition (\ref{EQ-355}) with $s=1-n$ that
\begin{eqnarray}
Y(1-n) &=&\int_{C}z^{-n-1}F_{q}(w,-z\mid w_{1})dz  \label{EQ-357} \\
&=&2\pi i\text{Res}_{z=0}\left\{ \text{ }z^{-n-1}F_{q}(w,-z\mid
w_{1})\right\}  \notag \\
&=&(2\pi i)\frac{(-1)^{n}}{n!}\beta _{n}(w:q\mid w_{1}),  \notag
\end{eqnarray}%
where we have made of the power-series representation in (\ref{EQ-3.1}).
Thus by Cauchy Residue Theorem, we easily \ complete the proof of Theorem 8
upon suitably combiningfin (\ref{EQ-358}) and (\ref{EQ-357}) with (\ref%
{EQ-356}).
\end{proof}

\begin{remark}
The representation in (\ref{EQ-3.4}) can be used to show that $\zeta
_{q}(s,w\mid w_{1})$ admits itself of an analytical continuation to whole
complex $s$- plane except for simple pole at $s=1$.
\end{remark}

\section{The Dirichlet Type Changhee $q$-$L$-Function}

We consider the following contour integral:
\begin{eqnarray}
\frac{1}{\Gamma (s)}\oint_{C}t^{s-2}F_{\chi ,q}(-t &\mid &w_{1})dt=\frac{1}{%
\Gamma (s)}\int_{0}^{\infty }t^{s-2}F_{\chi ,q}(-t\mid w_{1})dt  \notag \\
&=&\frac{w_{1}}{\Gamma (s)}\sum_{n=1}^{\infty }\chi
(n)q^{w_{1}n}\int_{0}^{\infty }t^{s-1}e^{-[w_{1}n]t}dt  \notag \\
&=&w_{1}\sum_{n=0}^{\infty }\frac{\chi (n)q^{w_{1}n}}{[w_{1}n]^{s}},
\label{EQ-4.1}
\end{eqnarray}%
where $C$ denote a positively oriented circle of radius $R$, centered at
origin. The function
\begin{equation*}
B(t)=\frac{1}{\Gamma (s)}t^{s-2}F_{\chi ,q}(-t\mid w_{1})
\end{equation*}
has pole $t=0$ inside the contour $C.$ Therefore, if want to integrate $B(t)$
function, than we have modify the contour by indentation at this point. We
take as indentation identical small semicircle, which has radius $r$,
leaving $t=0$.

We now define the Dirichlet's type Changhee $q$-$L$-function as follows:

\begin{definition}
Let $\chi $ be a Dirichlet character of conductor $f\in \mathbb{Z}^{+}$.
\begin{equation}
L_{q}(s,\chi \mid w_{1})=w_{1}\sum_{n=0}^{\infty }\frac{\chi (n)q^{w_{1}n}}{%
[w_{1}n]^{s}}.  \label{EQ-4.2}
\end{equation}
\end{definition}

We now give a relationship between $L_{q}(s,\chi \mid w_{1})$ and
generalized Changhee $q$-Bernoulli numbers. Thus we give the numbers $%
L_{q}(1-n,\chi \mid w_{1}),n\in \mathbb{Z}^{+}$, explicitly.

\begin{theorem}
Let $\chi $ be a Dirichlet character of conductor $f\in \mathbb{Z}^{+}$ and
let $n\in \mathbb{Z}^{+}$. Then
\begin{equation*}
L_{q}(1-n,\chi \mid w_{1})=-\frac{\beta _{n,\chi }(q\mid w_{1})}{n}.
\end{equation*}
\end{theorem}

\begin{proof}
Proof of Theorem 9 runs parallel to that of Theorem 8 above, so we choose to
omit the details involved.
\end{proof}

The Dirichlet's Type Changhee $q$-$L$-function and Hurwitz type Changhee $q$%
-zeta function are closely related, too. We gave this relation as follows:

\begin{theorem}
Let $\chi $ be a Dirichlet character of conductor $f\in \mathbb{Z}^{+}$.
Then
\begin{equation}
L_{q}(s,\chi \mid w_{1})=[f]^{-s}\sum_{a=1}^{f}\chi (a)\zeta _{q^{f}}(s,%
\frac{aw_{1}}{f}\mid w_{1}).  \label{EQ-4.3}
\end{equation}
\end{theorem}

\begin{proof}
By setting $n=a+kf$, where ( $k=0,1,2,...,\infty $ ; $a=1,2,...,f$ ) in (\ref%
{EQ-4.2}), we have
\begin{eqnarray*}
L_{q}(s,\chi &\mid &w_{1})=w_{1}\sum_{a=1}^{f}\chi (a)\sum_{k=0}^{\infty }%
\frac{q^{(aw_{1}+kfw_{1})}}{[aw_{1}+kfw_{1}]^{s}} \\
&=&w_{1}\sum_{a=1}^{f}\chi (a)\sum_{k=0}^{\infty }\frac{q^{f(\frac{aw_{1}}{f}%
+kw_{1})}}{[f]^{s}[\frac{aw_{1}}{f}+kw_{1}:q^{f}]^{s}} \\
&=&[f]^{-s}\sum_{a=1}^{f}\chi (a)\left\{ -\frac{(1-q^{f})^{s}}{s-1}\frac{1}{%
\log q^{f}}+w_{1}\sum_{n=0}^{\infty }\frac{q^{f(\frac{aw_{1}}{f}+kw_{1})}}{[%
\frac{aw_{1}}{f}+kw_{1}:q^{f}]^{s}}\right\} .
\end{eqnarray*}%
By using (\ref{EQ-3.5}) in the above, we easily arrive at the desired result
(\ref{EQ-4.3}).
\end{proof}

\section{Barnes' Type Changhee $q$-Bernoulli Numbers}

We now define new generating functions as follows:
\begin{eqnarray}
G_{q}(t &\mid &w_{1})=F_{q}(t\mid w_{1})-\frac{q-1}{\log q}\exp (\frac{t}{1-q%
})  \label{EQ-5.1} \\
&=&\sum_{k=0}^{\infty }B_{k}(q\mid w_{1})\frac{t^{k}}{k!}\text{ \ (}\mid
t\mid <2\pi \text{ ),}  \notag
\end{eqnarray}%
where the coefficients $B_{k}(q\mid w_{1})$ are called Barnes' type Changhee
$q$-Bernoulli numbers. By (\ref{EQ-5.1}), we easily see that
\begin{equation*}
B_{k}(q\mid w_{1})=\frac{q-1}{\log q}(\frac{1}{q-1})^{k}+\beta _{k}(q\mid
w_{1}),
\end{equation*}%
where $\beta _{k}(q\mid w_{1})$ is given by (\ref{EQ-3.1}).

Analogous to (\ref{EQ-2.1}), we can also consider the \textit{modified}
Changhee $q$-Bernoulli polynomials as follows:
\begin{eqnarray}
G_{q}(w,t &\mid &w_{1})=F_{q}(w,t\mid w_{1})-\frac{q-1}{\log q}\exp (\frac{t%
}{1-q})  \label{EQ-5.2} \\
&=&\sum_{n=0}^{\infty }\frac{B_{n}(w:q\mid w_{1})t^{n}}{n!}\text{ \ (}\mid
t\mid <2\pi \text{ ).}  \notag
\end{eqnarray}%
By applying (\ref{EQ-5.2}), we easily see that
\begin{eqnarray}
\frac{1}{\Gamma (s)}\int_{0}^{\infty }t^{s-2}G_{q}(w,-t &\mid
&w_{1})dt=w_{1}\sum_{n=0}^{\infty }\frac{q^{w_{1}n+w}}{[w_{1}n+w]^{s}}
\label{EQ-5.3} \\
&=&\zeta _{q,1}(s,w\mid w_{1})  \notag
\end{eqnarray}%
By using (\ref{EQ-5.2}), we give the following relationship between $\zeta
_{q,1}(s,w\mid w_{1})$ and $B_{n}(w:q\mid w_{1})$.

\begin{theorem}
For positive integer $n$,
\begin{equation*}
\zeta _{q,1}(1-n,w\mid w_{1})=-\frac{B_{n}(w:q\mid w_{1})}{n}.
\end{equation*}
\end{theorem}

We next define Barnes' type multiple Changhee $q$-Bernoulli polynomials as
follows:
\begin{eqnarray}
G_{q}^{(r)}(w,t &\mid &w_{1},w_{2},...,w_{r})  \notag \\
&=&(-t)^{r}(\prod_{i=1}^{r}w_{i})\sum_{n_{1},n_{2},...,n_{r}=0}^{\infty
}q^{w+n_{1}w_{1}+n_{2}w_{2}+...+n_{r}w_{r}}e^{[w+n_{1}w_{1}+n_{2}w_{2}+...+n_{r}w_{r}]t}
\notag \\
&=&\sum_{n=0}^{\infty }\frac{B_{n}^{(r)}(w:q\mid w_{1},w_{2},...,w_{r})t^{n}%
}{n!}\text{ \ (}\mid t\mid <2\pi \text{ ),}  \label{EQ-5.4}
\end{eqnarray}%
whith, as usual,
\begin{equation*}
\sum_{n_{1},n_{2},...,n_{r}=0}^{\infty }=\sum_{n_{1}=0}^{\infty
}\sum_{n_{2}=0}^{\infty }...\sum_{n_{r}=0}^{\infty }.
\end{equation*}%
It follows from (\ref{EQ-5.4}) that
\begin{equation*}
\lim_{q\rightarrow 1}G_{q}^{(r)}(w,t\mid w_{1},w_{2},...,w_{r})=\frac{%
e^{tw}(tw_{1})(tw_{2})...(tw_{r})}{%
(e^{tw_{1}}-1)(e^{tw_{2}}-1)...(e^{tw_{r}}-1)}.
\end{equation*}%
This gives generating function of Barnes' type multiple Bernoulli numbers.
Thus we get the following limit relationship:
\begin{equation*}
\lim_{q\rightarrow 1}B_{n}^{(r)}(w:q\mid
w_{1},w_{2},...,w_{r})=B_{n}^{(r)}(w\mid w_{1},w_{2},...,w_{r}).
\end{equation*}%
This gives Barnes' type multiple Bernoulli numbers as a limit when $q$
approaches.

By using (\ref{EQ-5.4}), we give Barnes' type Changhee multiple $q$-zeta
functions. For $s\in \mathbb{C}$, we consider the below integral which is
known Mellin transformation of $G_{q}^{(r)}(w,t\mid w_{1},w_{2},...,w_{r})$.
\begin{eqnarray}
\frac{1}{\Gamma (s)}\int_{0}^{\infty }t^{s-1-r}G_{q}^{(r)}(w,-t &\mid
&w_{1},w_{2},...,w_{r})dt  \notag \\
&=&(\prod_{i=1}^{r}w_{i})\sum_{n_{1},n_{2},...,n_{r}=0}^{\infty }\frac{%
q^{w+n_{1}w_{1}+n_{2}w_{2}+...+n_{r}w_{r}}}{%
^{[w+n_{1}w_{1}+n_{2}w_{2}+...+n_{r}w_{r}]^{s}}}.  \label{EQ-5.5}
\end{eqnarray}%
By using (\ref{EQ-5.5}), we can define Barnes' type Changhee multiple $q$%
-zeta functions as follows:

\begin{definition}
Let $s,w,w_{1},w_{2},...,w_{r}\in \mathbb{C}$ with $\func{Re}(w)>0$ and $%
r\in \mathbb{Z}^{+}$.
\begin{equation}
\zeta _{q,r}(s,w\mid
w_{1},w_{2},...,w_{r})=(\prod_{i=1}^{r}w_{i})%
\sum_{n_{1},n_{2},...,n_{r}=0}^{\infty }\frac{%
q^{w+n_{1}w_{1}+n_{2}w_{2}+...+n_{r}w_{r}}}{%
^{[w+n_{1}w_{1}+n_{2}w_{2}+...+n_{r}w_{r}]^{s}}}.  \label{EQ-5.7}
\end{equation}
\end{definition}

We note that $\zeta _{q,r}(s,w\mid w_{1},w_{2},...,w_{r})$ is analytic
continuation for $\func{Re}(s)>r$. By using (\ref{EQ-5.4}) and (\ref{EQ-5.5}%
), we arrive at the following theorem.

\begin{theorem}
Let $r\in \mathbb{Z}^{+}$. Then%
\begin{equation*}
\zeta _{q,r}(1-n,w\mid w_{1},w_{2},...,w_{r})=(-1)^{r}\frac{(n-1)!}{(n+r-1)!}%
B_{n+r-1}^{(r)}(w:q\mid w_{1},w_{2},...,w_{r}).
\end{equation*}
\end{theorem}

We record the following limit relationship:
\begin{eqnarray*}
\lim_{q\rightarrow 1}\zeta _{q,r}(1-n,w &\mid &w_{1},w_{2},...,w_{r})=\zeta
_{1,r}(1-n,w\mid w_{1},w_{2},...,w_{r}) \\
&=&(-1)^{r}\frac{(n-1)!}{(n+r-1)!}B_{n+r-1}^{(r)}(w\mid
w_{1},w_{2},...,w_{r})
\end{eqnarray*}%
between the ordinary Barnes' type multiple zeta functions and Barnes' type
Bernoulli numbers.

\section{Relations Between $\Gamma _{q}$, $\protect\zeta _{q,r}(s,w\mid
1,1,...,1)$ and $L_{q,r}(s,\protect\chi )$}

$\Gamma _{q}$-function is defined by ( see, for example, \cite{Kim8}, \cite%
{Kim13}, and \cite{Kim14})
\begin{equation}
\Gamma _{q}(z)=(1-q)^{1-z}\frac{(q:q)_{\infty }}{(q^{z}:q)_{\infty }},
\label{EQ-6.1}
\end{equation}%
where
\begin{equation*}
(q^{z}:q)_{\infty }=\prod_{k=0}^{\infty }(1-q^{k+z}).
\end{equation*}%
Thus the function $\zeta _{q}(s,x)$ is defined by
\begin{equation}
\zeta _{q}(s,x)=\sum_{n=0}^{\infty }\frac{q^{n+x}}{[n+x]^{s}}.
\label{EQ-6.2}
\end{equation}

If, for convenience, we denote
\begin{equation}
\frac{d}{ds}\zeta _{q}(s,x)\mid _{s=0}=\zeta _{q}^{^{\prime }}(0,x),
\label{EQ-6.3}
\end{equation}%
we readily observe that
\begin{equation}
\zeta _{q}(0,a+1)=\zeta _{q,1}^{^{\prime }}(0,1)+(q-1)\sum_{m=1}^{a}\log
[m]^{-s-1}-\sum_{n=0}^{\infty }[n]^{-s}  \label{EQ-6.4}
\end{equation}%
By using (\ref{EQ-6.1}), (\ref{EQ-6.2}) and (\ref{EQ-6.3}), we have the
following theorem.

\begin{theorem}
(\cite{Kim14})
\begin{equation}
\Gamma _{q}(a)=\frac{e^{\zeta _{q}^{^{\prime }}(0,a+1)}}{e^{\zeta
_{q}^{^{\prime }}(0,1)}}\prod_{m=1}^{a}[m]^{(1-q)[m]}.  \label{EQ-6.5}
\end{equation}
\end{theorem}

By means of (\ref{EQ-6.5}), we obtain
\begin{equation*}
\lim_{q\rightarrow 1}\Gamma _{q}(a)=\frac{e^{\zeta ^{^{\prime }}(0,a+1)}}{%
e^{\zeta ^{^{\prime }}(0,1)}}=\Gamma ^{^{\prime }}(a).
\end{equation*}

We now define $\zeta _{q,2}(s,x)$ as follows
\begin{equation}
\zeta _{q,2}(s,x)=\sum_{n_{1},n_{2}=0}^{\infty }\frac{q^{n_{1}+n_{2}+x}}{%
[n_{1}+n_{2}+x]^{s}}.  \label{EQ-6.6}
\end{equation}%
We give some properties of the zeta function defined by (\ref{EQ-6.6}) as
follows:

1)
\begin{equation}
\lbrack n]^{-s}\sum_{m=0}^{n-1}\sum_{k=0}^{n-1}\zeta _{q^{n},2}(s,x+\frac{k+m%
}{n})=\zeta _{q,2}(s,nx).  \label{EQ-6.7}
\end{equation}

2) By (\ref{EQ-6.7}), we have
\begin{equation}
\sum_{m=0}^{n-1}\sum_{k=0}^{n-1}\zeta _{q^{n},2}^{^{\prime }}(s,x+\frac{k+m}{%
n})=[n]^{s}\zeta _{q,2}(s,nx)\log [n]+[n]^{s}\zeta _{q,2}^{^{\prime }}(s,nx).
\label{EQ-6.8}
\end{equation}

3) In (\ref{EQ-6.8}), if we take $s=0$, so that we have
\begin{equation*}
\prod_{m=0}^{n-1}\prod_{k=0}^{n-1}e^{\zeta _{q^{n},2}^{^{\prime }}(0,x+\frac{%
k+m}{n})}=[n]^{\zeta _{q,2}(0,nx)}e^{\zeta _{q,2}^{^{\prime }}(0,nx)}.
\end{equation*}%
This function is called di-gamma function.

For any integer $k$ with $k\geq 0$, we define the function $W_{m,\chi ,q}(k)$
functions as follows:
\begin{equation*}
W_{m,\chi ,q}(k)=\sum_{a=1}^{k}\chi (a)q^{a}[a]^{m},\text{ }m\geq 0.
\end{equation*}%
If $\chi \equiv 1$ in the above, we have
\begin{equation*}
W_{m,q}(k)=W_{m,1,q}(k)=\sum_{a=1}^{k}q^{a}[a]^{m}.
\end{equation*}

By using (\ref{EQ-2.1}) and (\ref{EQ-2.6}), we easily obtain
\begin{equation*}
\frac{1}{t}(F_{q,\chi }(x,t)-F_{q,\chi }(t))=\sum_{k=0}^{\infty }(\frac{%
B_{k+1,\chi }(x:q)-B_{k+1,\chi }(q)}{k+1})\frac{t^{k}}{k!}.
\end{equation*}%
Thus we arrive at the following results:
\begin{equation*}
W_{k,\chi ,q}(n)=\frac{B_{k+1,\chi }(n:q)-B_{k+1,\chi }(q)}{k+1}
\end{equation*}%
and
\begin{equation*}
W_{k,q}(n)=\frac{B_{k+1}(n:q)-B_{k+1}(q)}{k+1}.
\end{equation*}

\section{Analytic Properties of $q$-$L$-Function and the $q$-Hurwitz zeta
function}

Let $s\in \mathbb{C}$. $q$-zeta function is defined by\cite{kim-Jang-Rim-Son}
\begin{equation}
\zeta _{q}(s)=-\frac{(1-q)^{s}}{s-1}\frac{1}{\log q}+\sum_{n=0}^{\infty }%
\frac{q^{n}}{[n]^{s}}.  \label{EQ-7.1}
\end{equation}%
We note that $\zeta _{q}(s)$ are analytically continued for $\func{Re}(s)>1$.

$q$-Hurwitz zeta function is defined by
\begin{equation}
\zeta _{q}(s,x)=-\frac{(1-q)^{s}}{s-1}\frac{1}{\log q}+\sum_{n=0}^{\infty }%
\frac{q^{n+x}}{[n+x]^{s}}.  \label{EQ-7.2}
\end{equation}%
From the definition (\ref{EQ-7.2}), one can easily get
\begin{equation*}
\zeta _{q}(0,v)=-(\frac{q-1}{\log q})[v]+\frac{q^{v}}{\log q}-\frac{^{q^{v}}%
}{q-1}.
\end{equation*}%
If $q\rightarrow 1$, then we have
\begin{equation}
\zeta _{q}(0,v)\rightarrow \zeta (0,v)=\frac{1}{2}-v.  \label{EQ-7.3}
\end{equation}%
More generally, we have
\begin{equation}
\zeta _{q}(-m,v)=-\frac{\beta _{m+1}(v:q)}{m+1},\text{ }m\geq 0.
\label{EQ-7.6}
\end{equation}

For $s\in \mathbb{C}$, we define $q$-$L$-function as follows:

Let $\chi $ be a Dirichlet character of conductor $f\in \mathbb{Z}^{+}$.
\begin{equation}
L_{q}(s,\chi )=\sum_{n=0}^{\infty }\frac{\chi (n)q^{n}}{[n]^{s}}.
\label{EQ-7.4}
\end{equation}

By using (\ref{EQ-7.1}) and (\ref{EQ-7.4}), we easily obtain the following
relation:
\begin{equation}
L_{q}(s,\chi )=[f]^{-s}\sum_{a=1}^{f}\chi (a)\zeta _{q^{f}}(s,\frac{a}{f}).
\label{EQ-7.5}
\end{equation}%
By using (\ref{EQ-7.5}) and (\ref{EQ-7.6}), we easily have the following
theorem.

\begin{theorem}
Let $k\in \mathbb{Z}^{+}$. We have
\begin{equation*}
L_{q}(1-k,\chi )=-\frac{\beta _{k,\chi }(q)}{k}.
\end{equation*}
\end{theorem}

In particular, if we define
\begin{equation*}
H_{q}(s,a,F)=\sum_{%
\begin{array}{c}
m>0 \\
m\equiv a(\func{mod}F)%
\end{array}%
}\frac{q^{m}}{[m]^{s}}.
\end{equation*}%
then we have
\begin{equation*}
H_{q}(s,a,F)=\sum_{n=0}^{\infty }\frac{q^{nF+a}}{[nF+a]^{s}}=\frac{1}{[F]^{s}%
}\zeta _{q^{F}}(s,\frac{a}{F}).
\end{equation*}%
It is well-known for the Hurwitz zeta function that (cf., e.g., \cite{H. M.
Srivastava and J. Choi}, p. 91, Equation (15)), we have
\begin{equation*}
\lim_{s\rightarrow \infty }(\zeta (s,a)-\frac{1}{s-1})=-\frac{\Gamma
^{^{\prime }}(a)}{\Gamma (a)}=-\psi (a)
\end{equation*}%
in terms of the familiar digamma ( or $\psi $) function. Hence, $%
H_{q}(s,a,F) $ has a simple pole at $s=1$ with residue
\begin{equation*}
\frac{1}{[F]}\frac{1}{F}\frac{q^{F}-1}{\log q}.
\end{equation*}

\begin{remark}
Barnes-Changhee multiple $q$-zeta functions are defined by (see\ \cite{Kim7}%
, \cite{Kim15})
\begin{equation*}
\zeta _{q,r}(s,w\mid
a_{1},a_{2},...,a_{r})=\sum_{n_{1},n_{2},...,n_{r}=0}^{\infty }\frac{%
q^{w+n_{1}+n_{2}+...+n_{r}}}{^{[w+n_{1}a_{1}+n_{2}a_{2}+...+n_{r}a_{r}]^{s}}}%
,\text{ }\Re (w)>0,\text{ }q\in C\text{ \ with }\mid q\mid <1\text{,}
\end{equation*}%
which, for $a_{1}=a_{2}=...=a_{r}=1,$ yields
\begin{equation*}
\zeta _{q,r}(s,w\mid 1,1,...,1)=\sum_{n_{1},n_{2},...,n_{r}=0}^{\infty }%
\frac{q^{w+n_{1}+n_{2}+...+n_{r}}}{^{[w+n_{1}+n_{2}+...+n_{r}]^{s}}}.
\end{equation*}%
Moreover, if $w=r$ and $\ s=1-n$ $\ $( $n\in \mathbb{Z}^{+}$ ), we have
\begin{equation*}
\zeta _{q,r}(1-n,r\mid 1,1,...,1)=(-1)^{r}\frac{(n-1)!}{(n+r-1)!}%
B_{n+r-1}^{(r)}(r:q).
\end{equation*}%
We also note that
\begin{eqnarray*}
\lim_{q\rightarrow 1}\zeta _{q,r}(1-n,r &\mid &1,1,...,1)=\zeta
_{r}(1-n,r\mid 1,1,...,1) \\
&=&(-1)^{r}\frac{(n-1)!}{(n+r-1)!}B_{n+r-1}^{(r)}(r).
\end{eqnarray*}%
(see\ \cite{Kim7}, \cite{Kim15}).

Similarly, by using the analogous approaches to the multiple $L$-functions,
we have
\begin{equation}
L_{q,r}(s,\chi )=\sum_{n_{1},n_{2},...,n_{r}=1}^{\infty }\frac{\chi
(n_{1})\chi (n_{2})...\chi (n_{r})q^{n_{1}+n_{2}+...+n_{r}}}{%
^{[n_{1}+n_{2}+...+n_{r}]^{s}}}.  \label{EQ-7.7}
\end{equation}%
For $\ s=-n$ $\ $( $n\in \mathbb{Z}^{+}$ ), we find from (\ref{EQ-3.3}) that
\begin{equation*}
L_{q,r}(-n,\chi )=(-1)^{r}\frac{n!}{(n+r)!}B_{n+r,\chi }^{(r)}(q)
\end{equation*}%
(see\ \cite{Kim7}, \cite{Kim15}).
\end{remark}

The following theorem provides a relationship between $L_{q,r}(s,\chi )$ and
$\zeta _{q,r}(s,w\mid a_{1},a_{2},...,a_{r})$.

\begin{theorem}
Let $\chi $ be a Dirichlet character of conductor $f\in \mathbb{Z}^{+}$.
Also let $a_{1},a_{2},...,a_{r}$ and $n_{1},...,n_{r}$ \ be in $\mathbb{Z}%
^{+}$. Then
\begin{equation*}
L_{q,r}(s,\chi )=[f]^{-s}\sum_{a_{1},...,a_{r}=1}^{f}\chi (a_{1})\chi
(a_{2})...\chi (a_{r})\zeta _{q^{f},r}(s,\frac{a_{1}+...+a_{r}}{f}\mid
1,...,1).
\end{equation*}
\end{theorem}

\begin{proof}
By setting $n_{j}=a_{j}+n_{j}f$, ( $j\in \left\{ 1,2,...,r\right\} $, $%
n_{j}=0,1,...,\infty $, and $a_{j}=1,2,...,f$ \ ) in (\ref{EQ-7.7}), we have
\begin{eqnarray*}
L_{q,r}(s,\chi ) &=&\sum_{a_{1},...,a_{r}=1}^{f}\chi (n_{1})\chi
(n_{2})...\chi (n_{r})\sum_{n_{1},n_{2},...,n_{r}=0}^{\infty }\frac{%
q^{a_{1}+...+a_{r}+f(n_{1}+n_{2}+...+n_{r})}}{%
^{[a_{1}+...+a_{r}+f(n_{1}+n_{2}+...+n_{r})]^{s}}} \\
&=&[f]^{-s}\sum_{a_{1},...,a_{r}=1}^{f}\chi (n_{1})\chi (n_{2})...\chi
(n_{r})\sum_{n_{1},n_{2},...,n_{r}=0}^{\infty }\frac{q^{f(\frac{%
a_{1}+...+a_{r}}{f}+n_{1}+n_{2}+...+n_{r})}}{^{[\frac{a_{1}+...+a_{r}}{f}%
+n_{1}+n_{2}+...+n_{r}]^{s}}} \\
&=&[f]^{-s}\sum_{a_{1},...,a_{r}=1}^{f}\chi (n_{1})\chi (n_{2})...\chi
(n_{r})\zeta _{q^{f},r}(s,\frac{a_{1}+...+a_{r}}{f}\mid 1,1,...,1).
\end{eqnarray*}%
Thus we obtain the desired result asserted by Theorem 15.
\end{proof}

By putting $w_{1}=w_{2}=...=w_{r}=1$ and $w=x$ in (\ref{EQ-5.7}), we obtain
\begin{equation}
\zeta _{q,r}(s,x\mid 1,1,...,1)=\sum_{l=0}^{\infty }\left(
\begin{array}{c}
l+r-1 \\
r-1%
\end{array}%
\right) \frac{q^{x+l}}{[x+l]^{s}},  \label{EQ-7.8}
\end{equation}%
which, for $r=1$, reduces immediately to the $q$-Hurwitz zeta function:
\begin{equation*}
\zeta _{q,r}(s,x\mid 1,1,...,1)=\sum_{l=0}^{\infty }\frac{q^{x+l}}{[x+l]^{s}}%
.
\end{equation*}%
Furthermore, by setting $s=-n$, with ( $n\in \mathbb{Z}^{+}$ ) in (\ref%
{EQ-7.8}), we have
\begin{eqnarray*}
\zeta _{q,r}(-n,x &\mid &1,1,...,1)=\sum_{l=0}^{\infty }\left(
\begin{array}{c}
l+r-1 \\
r-1%
\end{array}%
\right) \frac{q^{x+l}}{[x+l]^{-n}} \\
&=&(-1)^{r}\frac{(n-1)!}{(n+r-1)!}\beta _{n+r-1}^{(r)}(x:q\mid 1,1,...,1),
\end{eqnarray*}%
where $B_{n}^{(r)}(x:q\mid 1,1,...,1)$ numbers are defined as follows (\cite%
{Kim10}, \cite{Kim11}):

For $n,k\in \mathbb{Z}^{+}$ ($k>1$), if $S_{n,q}$ denotes the sums of the $n$%
th powers of positive $q$-integers up to $k-1$ (see\cite{Kim11}):
\begin{equation*}
S_{n,q}(k)=\sum_{m=0}^{k-1}q^{m}[m]^{n},
\end{equation*}%
then we have
\begin{equation*}
B_{n,q}^{(r)}(x:q\mid 1,1,...,1)=(-1)^{r}\frac{(r+1)!}{r!}\sum_{k=0}^{\infty
}\left(
\begin{array}{c}
k+r-1 \\
r-1%
\end{array}%
\right) q^{x+k}\sum_{m=0}^{r-1}S_{m,q^{r-m}}(x+k).
\end{equation*}

\begin{remark}
By using (\ref{EQ-5.7}), we have
\begin{equation*}
\lim_{q\rightarrow \infty }\zeta _{q,r}(s,w\mid a_{1},a_{2},...,a_{r})=\zeta
_{r}(s,w\mid a_{1},a_{2},...,a_{r})
\end{equation*}%
in termes of Barnes multiple zeta function in (\ref{Eq-12}).
\end{remark}

\section{Generalized Multiple Changhee $q$-Bernoulli Numbers and the
Dirichlet Type Multiple Changhee $q$-$L$-functions}

In this section, we define generalized multiple Changhee $q$-Bernoulli
numbers attached to the Drichlet character $\chi $. We also construct
Dirichlet's type multiple Changhee $q$-$L$-functions. We then give relation
between Dirichlet's type multiple Changhee $q$-$L$-functions and generalized
multiple Changhee $q$-Bernoulli numbers as well.

The generalized multiple Changhee $q$-Bernoulli numbers attached to the
Drichlet character $\chi $ are defined by means of the following generating
function:
\begin{eqnarray}
F_{q,\chi }^{(r)}(t &\mid &w_{1},...,w_{r})=(-t)^{r}\left(
\prod_{j=1}^{r}w_{j}\right) \sum_{n_{1},...,n_{r}=1}^{\infty }\left(
\prod_{k=1}^{r}\chi (n_{k})\right) q^{\left( \sum_{m=1}^{r}w_{m}n_{m}\right)
}e^{[\sum_{m=1}^{r}w_{m}n_{m}]t}  \notag \\
&=&\sum_{n=0}^{\infty }B_{n,\chi }^{(r)}(q\mid w_{1},...,w_{r})\frac{t^{n}}{%
n!}\text{ \ (}\mid t\mid <2\pi \text{ ),}  \label{EQ-8.1}
\end{eqnarray}%
where $w_{1},...,w_{r}\in \mathbb{R}^{+}\mathbb{,}$ $r\in \mathbb{Z}^{+}$.

By simple calculations in (\ref{EQ-8.1}), we have
\begin{eqnarray}
&&(-t)^{r}\left( \prod_{j=1}^{r}w_{j}\right)
\sum_{n_{1},...,n_{r}=1}^{\infty }\left( \prod_{k=1}^{r}\chi (n_{k})\right)
q^{\left( \sum_{m=1}^{r}w_{m}n_{m}\right) }e^{[\sum_{m=1}^{r}w_{m}n_{m}]t}
\notag \\
&=&(-t)^{r}(\prod_{j=1}^{r}w_{j})\sum_{a_{1},...,a_{r}=1}^{f}\left(
\prod_{k=1}^{r}\chi (n_{k})q^{a_{k}w_{k}}\right) \exp
(t[\sum_{m=1}^{r}w_{m}a_{m}])  \notag \\
&&.\sum_{n_{1},...,n_{r}=1}^{\infty }q^{\left(
\sum_{m=1}^{r}w_{m}n_{m}f\right) }\exp \left(
t[\sum_{m=1}^{r}w_{m}n_{m}f]q^{\sum_{m=1}^{r}w_{m}n_{m}f}\right)  \notag \\
&=&F_{q,\chi }^{(r)}(t\mid w_{1},...,w_{r})  \label{EQ-8.2}
\end{eqnarray}%
Now by applying (\ref{EQ-5.4}), we obtain
\begin{equation}
F_{q,\chi }^{(r)}(t\mid
w_{1},...,w_{r})=[f]^{-r}\sum_{a_{1},...,a_{r}=1}^{f}\prod_{i=1}^{r}\chi
(a_{i})G_{q^{f}}^{(r)}(\frac{w_{1}a_{1}+...+w_{r}a_{r}}{f},[f]t\mid
w_{1},w_{2},...,w_{r}).  \label{EQ-8.3}
\end{equation}

By using (\ref{EQ-5.4}) and (\ref{EQ-8.3}), we readily arrive at the
following theorem:

\begin{theorem}
Let $\chi $ be a Dirichlet character of conductor $f\in \mathbb{Z}^{+}$.
Then
\begin{equation}
B_{n,\chi }^{(r)}(q\mid
w_{1},...,w_{r})=[f]^{n-r}\sum_{a_{1},...,a_{r}=1}^{f}\prod_{i=1}^{r}\chi
(a_{i})B_{n}^{(r)}(\frac{w_{1}a_{1}+...+w_{r}a_{r}}{f}:q^{f}\mid
w_{1},w_{2},...,w_{r}).  \label{EQ-8.4}
\end{equation}
\end{theorem}

Here, we can now construct Dirichlet's type multiple Changhee $q$-$L$%
-function. By using Mellin transformation and Residue Theorem in (\ref%
{EQ-8.1}), then we obtain
\begin{eqnarray}
\frac{1}{\Gamma (s)}\int_{0}^{\infty }t^{s-1-r}F_{q,\chi }^{(r)}(-t &\mid
&w_{1},...,w_{r})dt  \notag \\
&=&\left( \prod_{j=1}^{r}w_{j}\right) \sum_{n_{1},...,n_{r}=1}^{\infty
}\left( \prod_{k=1}^{r}\chi (n_{k})\right) q^{\left(
\sum_{m=1}^{r}w_{m}n_{m}\right) }\frac{1}{\Gamma (s)}\int_{0}^{\infty
}t^{s-1}e^{-[\sum_{m=1}^{r}w_{m}n_{m}]t}dt  \label{EQ-8.5}
\end{eqnarray}%
By using (\ref{EQ-8.5}), we can define Dirichlet's type multiple Changhee $q$%
-$L$-functions as follows.

\begin{definition}
For a Dirichlet character $\chi $with conductor $f\in \mathbb{Z}^{+}$, we
define
\begin{equation}
L_{q,r}(s,\chi \mid w_{1},...,w_{r})=\left( \prod_{j=1}^{r}w_{j}\right)
\sum_{n_{1},n_{2},...,n_{r}=1}^{\infty }\frac{\left( \prod_{k=1}^{r}\chi
(n_{k})\right) q^{\left( \sum_{m=1}^{r}w_{m}n_{m}\right) }}{%
^{[\sum_{m=1}^{r}w_{m}n_{m}]^{s}}}.  \label{EQ-8.6}
\end{equation}
\end{definition}

By (\ref{EQ-8.5}) and (\ref{EQ-8.6}), we can easily obtain relationship the
following relationship between $L_{q,r}(s,\chi \mid w_{1},...,w_{r})$ and $%
\zeta _{r,q}(s,w_{1}a_{1}+...+w_{r}a_{r}\mid w_{1},...,w_{r})$.

\begin{theorem}
\begin{equation}
L_{q,r}(s,\chi \mid
w_{1},...,w_{r})=[f]^{r-s}\sum_{a_{1},...,a_{r}=1}^{f}\left(
\prod_{k=1}^{r}\chi (a_{k})\right) \zeta _{q^{f},r}(s,\frac{%
w_{1}a_{1}+...+w_{r}a_{r}}{f}\mid w_{1},...,w_{r}).  \label{EQ-8.7}
\end{equation}
\end{theorem}

By using (\ref{EQ-8.1}) to (\ref{EQ-8.7}), the numbers $L_{q,r}(-n,\chi \mid
w_{1},...,w_{r})$, ( $n>0$ ) are given explicitly aby Theorem 18 bellow.

\begin{theorem}
\begin{equation*}
L_{q,r}(-n,\chi \mid w_{1},...,w_{r})=(-1)^{r}\frac{n!}{(n+r)!}B_{n,\chi
}^{(r)}(q\mid w_{1},...,w_{r}).
\end{equation*}
\end{theorem}

\section{Euler-Barnes' type Multiple $q$-Daehee Zeta Functions}

The main purpose of this section is to prove analytic continuation of the
Euler-Barnes' type multiple $q$-Daehee zeta functions depending on the
parameters $a_{1},...,a_{r}$ which are taken positive parts in the Complex
Field. Therefore, we construct generating function of $q$-Euler-Barnes'type
multiple Frobenius-Euler polynomials. We define Euler-Barnes' type multiple $%
q$-Daehee zeta Functions. As we remarked earlier, the Euler-Barnes' type
multiple $q$-Daehee zeta functions have a potentially useful connection with
Topology and Physics, together with the algebraic relations among them. We
give the values of these functions at negative integers as well.

Recently, by using an invariant $p$-adic integral, Kim\cite{Kim4} showed
that the Daehee numbers are related to the $q$-Bernoulli and Eulerian
numbers. Here we construct our generating function in complex case.

We now define the generating function of Euler-Barnes' type $q$-Daehee
numbers as follows.

If $w_{1}\in \mathbb{C}$ with positive real part and $u\in \mathbb{C}$ with $%
\mid u\mid <1$, then
\begin{eqnarray}
F_{u^{-1},q}(t &\mid &w_{1})=(1-u)\exp \left( \frac{t}{1-q}\right)
\sum_{j=0}^{\infty }\left( \frac{1}{1-q}\right) ^{j}\left( \frac{1}{%
1-q^{w_{1}j}u}\right) \frac{(-t)^{j}}{j!}  \notag \\
&=&\sum_{n=0}^{\infty }H_{n}(u^{-1}:q\mid w_{1})\frac{t^{n}}{n!}\text{ \ (}%
\mid t\mid <2\pi \text{ ).}  \label{EQ-9.1}
\end{eqnarray}%
By letting $q\rightarrow 1$ in (\ref{EQ-9.1}), we arrive at (\ref{eq-5}),
that is,
\begin{equation*}
\lim_{q\rightarrow 1}F_{u^{-1},q}(t\mid w_{1})=\frac{1-u^{-1}}{e^{wt}-u^{-1}}%
=\sum_{n=0}^{\infty }H_{n}(u^{-1}\mid w_{1})\frac{t^{n}}{n!},
\end{equation*}%
which implies that
\begin{equation*}
\lim_{q\rightarrow 1}H_{n}(u^{-1}:q\mid w_{1})=H_{n}(u^{-1}\mid w_{1}).
\end{equation*}%
By using (\ref{EQ-9.1}), we also get
\begin{equation}
F_{u^{-1},q}(t\mid w_{1})=(1-u)\sum_{k=0}^{\infty }\sum_{n=0}^{k}\left(
\begin{array}{c}
k \\
n%
\end{array}%
\right) \frac{(-1)^{n}}{1-q^{w_{1}n}u}\frac{1}{(1-q)^{k}}\frac{t^{k}}{k!}.
\label{EQ-9.2}
\end{equation}%
By applying (\ref{EQ-9.1}) and (\ref{EQ-9.2}), we easily obtain the
following result:

For $w_{1}\in \mathbb{C}$ with positive real part, $u\in \mathbb{C}$ with $%
\mid u\mid <1$ and $\func{Re}(w_{1})>0$,%
\begin{equation}
H_{k}(u^{-1}:q\mid w_{1})=\frac{(1-u)}{(1-q)^{k}}\sum_{n=0}^{k}\left(
\begin{array}{c}
k \\
n%
\end{array}%
\right) \frac{(-1)^{n}}{1-q^{w_{1}n}u}.  \label{EQ-9.3}
\end{equation}%
We note that the numbers in (\ref{EQ-9.3}) are called Euler-Barnes' type
Daehee $q$-Euler numbers.

By (\ref{EQ-9.1}), we have
\begin{eqnarray}
F_{u^{-1},q}(t &\mid &w_{1})=(1-u)e^{\frac{t}{1-q}}\sum_{j=0}^{\infty }(-%
\frac{1}{1-q})^{j}(\sum_{n=0}^{\infty }q^{w_{1}jn}u^{n})\frac{t^{j}}{j!}
\notag \\
&=&(1-u)\sum_{n=0}^{\infty }u^{n}e^{[w_{1}n]t}\text{ \ (}\mid t\mid <2\pi
\text{ ).}  \label{EQ-9.4}
\end{eqnarray}%
Using (\ref{EQ-9.1}) and (\ref{EQ-9.4}), we define generating function of
Euler-Barnes' type Daehee $q$-Euler polynomials as follows:
\begin{eqnarray}
F_{u^{-1},q}(t,w &\mid &w_{1})=e^{[w]t}F_{u^{-1},q}(q^{w}t\mid w_{1})
\label{EQ-9.5} \\
&=&\sum_{n=0}^{\infty }H_{n}(u^{-1},w:q\mid w_{1})\frac{t^{n}}{n!}\text{ \ (}%
\mid t\mid <2\pi \text{ ),}  \notag
\end{eqnarray}%
which implies that
\begin{equation}
F_{u^{-1},q}(t,w\mid w_{1})=(1-u)\sum_{n=0}^{\infty }u^{n}e^{[w+w_{1}n]t}.
\label{EQ-9.6}
\end{equation}%
Next we note that
\begin{eqnarray*}
\lim_{_{q\rightarrow \infty }}F_{u^{-1},q}(t,w &\mid
&w_{1})=(1-u)\sum_{n=0}^{\infty }u^{n}e^{(w+w_{1}n)t} \\
&=&\frac{(1-u^{-1})}{e^{w_{1}t}-u^{-1}}e^{wt} \\
&=&\sum_{n=0}^{\infty }H_{n}(u^{-1},w\mid w_{1})\frac{t^{n}}{n!}\text{ \ (}%
\mid t\mid <2\pi \text{ ),}
\end{eqnarray*}%
which gives (\ref{eq-6}). Hence
\begin{equation*}
\lim_{_{q\rightarrow \infty }}H_{n}(u^{-1},w:q\mid w_{1})=H_{n}(u^{-1},w\mid
w_{1}).
\end{equation*}

Now by applying (\ref{EQ-9.1}), (\ref{EQ-9.3}) , (\ref{EQ-9.5}) and (\ref%
{EQ-9.6}), we easily arrive at the following theorem.

\begin{theorem}
\begin{equation}
H_{n}(u^{-1},w:q\mid w_{1})=\sum_{l=0}^{n}\left(
\begin{array}{c}
n \\
l%
\end{array}%
\right) [w]^{n-l}q^{wl}H_{l}(u^{-1}:q\mid w_{1}).  \label{EQ-9.7}
\end{equation}
\end{theorem}

We remark that the numbers $H_{n}(u^{-1},w:q\mid w_{1})$ are called
Euler-Barnes' type Daehee $q$-Euler polynomials.

Let us define the Euler-Barnes' type multiple Daehee $q$-Euler polynomials.
The generating function of this polynomials are defined as follows:
\begin{eqnarray}
F_{u^{-1},q}^{(r)}(t,w &\mid
&w_{1},...,w_{r})=(1-u)^{r}\sum_{n_{1},...,n_{r}=0}^{\infty
}u^{n_{1}+...+n_{r}}e^{[w+w_{1}n_{1}+...+w_{r}n_{r}]t}  \notag \\
&=&\sum_{n=0}^{\infty }H_{u^{-1},q}^{(r)}(u^{-1},w\mid w_{1},...,w_{r})\frac{%
t^{n}}{n!}\text{ \ (}\mid t\mid <2\pi \text{ ).}  \label{EQ-9.8}
\end{eqnarray}%
where $r\in \mathbb{Z}^{+},$ $w_{1},...,w_{r}\in \mathbb{C}$ with positive
real part, $u\in \mathbb{C}$ with $\mid u^{-1}\mid <1.$

Note that that cf \cite{Kim9}
\begin{eqnarray}
\lim_{_{q\rightarrow 1}}F_{u^{-1},q}^{(r)}(t,w &\mid
&w_{1},...,w_{r})=(1-u)^{r}\sum_{n_{1},...,n_{r}=0}^{\infty
}u^{n_{1}+...+n_{r}}e^{(w+w_{1}n_{1}+...+w_{r}n_{r})t}  \notag \\
&=&\frac{(1-u^{-1})...(1-u^{-1})}{(e^{w_{1}t}-u^{-1})...(e^{w_{r}t}-u^{-1})}%
e^{wt}  \label{EQ-9.9} \\
&=&\sum_{n=0}^{\infty }H_{n}^{(r)}(u^{-1},w\mid w_{1},...,w_{r})\frac{t^{n}}{%
n!}\text{ \ (}\mid t\mid <2\pi \text{ ).}  \notag
\end{eqnarray}

By using (\ref{EQ-9.7}) to (\ref{EQ-9.9}), we have \cite{Kim9}
\begin{equation*}
\lim_{_{q\rightarrow \infty }}H_{n}^{(r)}(u^{-1},w:q\mid
w_{1},...,w_{r})=H_{n}^{(r)}(u^{-1},w\mid w_{1},...,w_{r}).
\end{equation*}

\section{Euler-Barnes' Type Daehee $q$-Zeta Functions}

By applying Mellin transformation and Residue Theorem in (\ref{EQ-9.1}) and (%
\ref{EQ-9.6}), we have
\begin{eqnarray}
\frac{1}{\Gamma (s)}\int_{0}^{\infty }t^{s-1}\frac{1}{1-u}F_{u^{-1},q}(-t,w
&\mid &w_{1})dt  \notag \\
&=&\sum_{n=0}^{\infty }u^{n}\frac{1}{\Gamma (s)}\int_{0}^{\infty
}e^{-[w+w_{1}n]t}t^{s-1}dt  \label{EQ-10.1} \\
&=&\sum_{n=0}^{\infty }\frac{u^{n}}{[w+w_{1}n]^{s}},  \notag
\end{eqnarray}%
where $\Gamma (s)$ is the Euler gamma function.

Thus by virtue of (\ref{EQ-10.1}), we consider Euler-Barnes' type Daehee $q$%
-zeta functions as follows.

For $s\in \mathbb{C}$,
\begin{equation}
\zeta _{q}(s,w,u\mid w_{1})=\sum_{n=0}^{\infty }\frac{u^{n}}{[w+w_{1}n]^{s}}.
\label{EQ-10.2}
\end{equation}%
We note that $\zeta _{q}(s,w,u\mid w_{1})$ is analytic for $\func{Re}(s)>1$,
and\ that $\zeta _{q}(s,u\mid w_{1})$ is called the Euler-Barnes' type
Daehee $q$-zeta functions which are defined as follows:
\begin{equation*}
\zeta _{q}(s,u\mid w_{1})=\sum_{n=0}^{\infty }\frac{u^{n}}{[w_{1}n]^{s}}.
\end{equation*}%
By using (\ref{EQ-10.1}) and (\ref{EQ-10.2}), we easily see that

\begin{theorem}
Let $n\in \mathbb{Z}^{+}$. Then
\begin{equation*}
\zeta _{q}(-n,w,u\mid w_{1})=\frac{1}{1-u}H_{n}(u^{-1},w:q\mid w_{1}).
\end{equation*}
\end{theorem}

By using the same method as in (\ref{EQ-10.1}), we shall construct the
analytic Euler-Barnes' type multiple Daehee $q$-zeta functions as follows:
\begin{eqnarray}
\frac{1}{\Gamma (s)}\int_{0}^{\infty }t^{s-1-r}\frac{1}{(1-u)^{r}}%
F_{u^{-1},q}^{(r)}(-t,w &\mid &w_{1},...,w_{r})dt  \notag \\
&=&\sum_{n_{1},n_{2},...,n_{r}=0}^{\infty }\frac{u^{n_{1}+n_{2}+...+n_{r}}}{%
^{[w+n_{1}w_{1}+n_{2}w_{2}+...+n_{r}w_{r}]^{s}}},  \label{EQ-10.3}
\end{eqnarray}%
where $\func{Re}(w)>0$ and $\func{Re}(s)>r$.

By means of (\ref{EQ-10.3}), we define Euler-Barnes' type multiple Daehee $q$%
-zeta functions.

\begin{definition}
\begin{equation*}
\zeta _{q}(s,w,u\mid
w_{1},w_{2},...,w_{r})=\sum_{n_{1},n_{2},...,n_{r}=0}^{\infty }\frac{%
u^{n_{1}+n_{2}+...+n_{r}}}{^{[w+n_{1}w_{1}+n_{2}w_{2}+...+n_{r}w_{r}]^{s}}},
\end{equation*}%
where $\Re (w)>0$ and $\Re (s)>r$.
\end{definition}

We note that for $\func{Re}(s)>1$, $\zeta _{q}(s,w,u\mid
w_{1},w_{2},...,w_{r})$ provides an analytic continuation in the Complex
Field $\mathbb{C}$ and that \cite{Kim9}
\begin{equation*}
\lim_{q\rightarrow 1}\zeta _{q}(s,w,u\mid w_{1},w_{2},...,w_{r})=\zeta
(s,w,u\mid w_{1},w_{2},...,w_{r}).
\end{equation*}

Analytic continuation and special values of Euler-Barnes' type multiple
Daehee $q$-zeta functions are given by integral representation of $\zeta
_{q}(s,w,u\mid w_{1},w_{2},...,w_{r})$ as follows:
\begin{equation*}
\zeta _{q}(s,w,u\mid w_{1},w_{2},...,w_{r})=\frac{1}{\Gamma (s)}%
\int_{0}^{\infty }t^{s-1-r}\frac{1}{(1-u)^{r}}F_{u^{-1},q}^{(r)}(-t,w\mid
w_{1},...,w_{r})dt.
\end{equation*}%
Now, by using Cauchy Residue Theorem for $s=n$ ( $n\in \mathbb{Z}^{+}$ ), we
obtain the following theorem.

\begin{theorem}
Let $n\in \mathbb{Z}^{+}$. Then%
\begin{equation*}
\zeta _{q}(-n,w,u\mid w_{1},w_{2},...,w_{r})=\frac{1}{(1-u)^{r}}%
H_{n}^{(r)}(u^{-1},w:q\mid w_{1},...,w_{r}).
\end{equation*}
\end{theorem}

\section{Further Remarks and Observations}

Shiratani\cite{K. Shiratani} defined the following zeta functions:
\begin{equation*}
\zeta (s\mid u)=\sum_{n=0}^{\infty }\frac{u^{-n}}{n^{s}},\text{ ( }\Re (s)>1%
\text{ ).}
\end{equation*}%
By using (\ref{eq-5}),the values of this function at negative integers are
obtained explicitly as follows:
\begin{equation*}
\zeta (-k\mid u)=-\frac{H_{k}(u)}{k},\text{ ( }k\in \mathbb{Z}^{+}\text{ ).}
\end{equation*}%
This functions generalized by Kim \cite{Kim9} to the form which is given
already in (\ref{eqq-20}). He also gave the analytic continuation of
multiple zeta functions ( the Euler-Barnes multiple zeta functions )
depending on parameters $a_{1},...,a_{r}$ taking positive values in the
complex number field. If \ $q\rightarrow 1$ in (\ref{EQ-9.9}), we get the $r$%
th Frobenius -Euler polynomials with parameters $w,w_{1},...,w_{r}$ taking
positive values in the complex number field
\begin{equation*}
\frac{(1-u^{-1})...(1-u^{-1})}{(e^{w_{1}t}-u^{-1})...(e^{w_{r}t}-u^{-1})}%
e^{wt}=\sum_{n=0}^{\infty }H_{n}^{(r)}(u^{-1},w\mid w_{1},...,w_{r})\frac{%
t^{n}}{n!}\text{ \ (}\mid t\mid <2\pi \text{ ),}
\end{equation*}%
If we set $w=0$ in the above generating function, then we readily obtain
multiple Euler-Barnes numbers given by
\begin{equation*}
H_{n}^{(r)}(u^{-1},0\mid w_{1},...,w_{r})=H_{n}^{(r)}(u^{-1}\mid
w_{1},...,w_{r}).
\end{equation*}%
We now define
\begin{eqnarray}
F_{u,q}^{(r)}(t,x\mid w_{1},...,w_{r})= &&\frac{(1-u)^{r}e^{-xt}}{%
(e^{-w_{1}t}-u)...(e^{-w_{r}t}-u)}(\frac{u}{u-1})^{r}  \notag \\
&=&(\frac{u^{f}}{u^{f}-1})^{r}\sum_{a_{1},...,a_{r}=1}^{f-1}u^{-%
\sum_{j=1}^{r}a_{j}}\frac{(1-u^{f})^{r}\exp (-ft\frac{x+%
\sum_{j=1}^{r}a_{j}w_{j}}{f})}{(e^{-w_{1}ft}-u^{f})...(e^{-w_{r}ft}-u^{f})}
\label{Rem-1}
\end{eqnarray}%
By using (\ref{Rem-1}), we have
\begin{eqnarray*}
\frac{u^{r}}{(u-1)^{r}}H_{n}^{(r)}(u,x &\mid &w_{1},...,w_{r}) \\
&=&f^{n}\sum_{a_{1},...,a_{r}=1}^{f}\frac{u^{fr-\sum_{j=1}^{r}a_{j}}}{%
(u^{f}-1)^{r}}H_{n}^{(r)}(u^{f},\frac{x+a_{1}w_{1}+...+a_{r}w_{r}}{f}\mid
w_{1},...,w_{r}).
\end{eqnarray*}%
This is known as distribution function. Now, by using this function,
Frobenius-Barnes' type measure is defined as follows.

Let $\chi $ be a Dirichlet character with conductor $f\in \mathbb{Z}^{+}$.
We define
\begin{eqnarray*}
&&\sum_{a_{1},...,a_{r}=1}^{f}\prod_{j=1}^{r}\chi
(a_{j})u^{rf-\sum_{j=1}^{r}a_{j}}\frac{e^{t%
\sum_{j=1}^{r}a_{j}w_{j}}(1-u^{f})^{r}}{%
(e^{w_{1}ft}-u^{f})...(e^{w_{r}ft}-u^{f})} \\
&=&\sum_{n=0}^{\infty }H_{n,\chi }^{(r)}(u\mid w_{1},...,w_{r})\frac{t^{n}}{%
n!} \\
&=&\sum_{n=0}^{\infty }(f^{n}\sum_{a_{1},...,a_{r}=1}^{f}\prod_{j=1}^{r}\chi
(a_{j})u^{rf-\sum_{j=1}^{r}a_{j}}H_{n}^{(r)}(u^{f},\frac{%
a_{1}w_{1}+...+a_{r}w_{r}}{f}\mid w_{1},...,w_{r}))\frac{t^{n}}{n!}.
\end{eqnarray*}%
By comparing coefficients $\frac{t^{n}}{n!}$, we easily see that
\begin{eqnarray*}
H_{n,\chi }^{(r)}(u &\mid &w_{1},...,w_{r}) \\
&=&f^{n}\sum_{a_{1},...,a_{r}=1}^{f}\prod_{j=1}^{r}\chi
(a_{j})u^{rf-\sum_{j=1}^{r}a_{j}}H_{n}^{(r)}(u^{f},\frac{%
a_{1}w_{1}+...+a_{r}w_{r}}{f}\mid w_{1},...,w_{r}).
\end{eqnarray*}%
By using this generating function, we\ can obtain an analytic continuation
of $\zeta (s,w,u\mid w_{1},w_{2},...,w_{r})$.

By using Mellin transformation in (\ref{Rem-1}), we easily see that
\begin{equation*}
\zeta _{r}(s,w,u\mid w_{1},w_{2},...,w_{r})=\frac{1}{\Gamma (s)}%
\int_{0}^{\infty }t^{s-1-r}\frac{1}{(1-u)^{r}}F_{u,q}^{(r)}(-t,w\mid
w_{1},...,w_{r})dt.
\end{equation*}%
Putting $s=-n$ ( $n>0$ ) in the above, we have
\begin{equation*}
\zeta _{r}(-n,w,u\mid w_{1},w_{2},...,w_{r})=\frac{u^{r}}{(u-1)^{r}}%
H_{n}^{(r)}(u,w\mid w_{1},...,w_{r}).
\end{equation*}

In $p$-adic case, we similarly obtain the following results.

For $u\in \mathbb{C}_{p}$, with $\mid 1-u\mid _{p}\geq 1$ we consider the
integral
\begin{equation*}
\lim_{N\rightarrow \infty }\frac{1}{1-u^{p^{N}}}%
\sum_{j=0}^{p^{N}-1}u^{p^{N}-j}g(j)=\sum_{j=0}^{p^{N}-1}g(j)E_{u}(j+p^{N}%
\mathbb{Z}_{p}).
\end{equation*}%
For $g\in UD(\mathbb{Z}_{p},\mathbb{C}_{p})$, the above limit exist. Thus
Euler integral defined as follows:
\begin{equation*}
\int_{\mathbb{Z}_{p}}g(x)dE_{u}(x)=\lim_{N\rightarrow \infty
}\sum_{x=0}^{p^{N}-1}g(x)\frac{u^{p^{N}-x}}{1-u^{p^{N}}}.
\end{equation*}

We now define $p$-adic Barnes' type Frobenius-Euler measure as follows:

Let $w_{1},w_{2},...,w_{r}$ be the nonzero $p$-adic integers. Thus we have
\begin{equation*}
E_{u,w_{1}}^{(k)}(x+p^{N}\mathbb{Z}_{p})=\frac{u^{p^{N}-x}}{1-u^{p^{N}}}%
H_{k}^{(1)}(u^{p^{N}},\frac{w_{1}x}{p^{N}}\mid w_{1}).
\end{equation*}%
This is a measure, because it is easily observe that
\begin{equation*}
\sum_{i=0}^{p^{N}-1}E_{u,w_{1}}^{(k)}(x+ip^{N}+p^{N+1}\mathbb{Z}%
_{p})=E_{u,w_{1}}^{(k)}(x+p^{N}\mathbb{Z}_{p}).
\end{equation*}%
Thus $E_{u,w_{1}}^{(k)}$ is a distribution. Now we give bounded property of%
\begin{equation*}
E_{u,w_{1}}^{(k)}(x+p^{N}\mathbb{Z}_{p}),\text{ when }\mid 1-u\mid _{p}\geq
1.
\end{equation*}%
Hence $E_{u,w_{1}}^{(k)}$ is a measure on $\mathbb{Z}_{p}$. By using this
measure and the above relations, we have
\begin{equation}
\int_{\mathbb{X}}\chi (x)dE_{u,w_{1}}^{(k)}(x)=\frac{1}{1-u^{f}}H_{k,\chi
}^{(1)}(u\mid w_{1}).  \label{Rem-2}
\end{equation}%
By using simple calculation, we see that
\begin{equation*}
\int_{\mathbb{X}}dE_{u,w_{1}}^{(k)}(x)=w_{1}^{k}\int_{\mathbb{X}}dE_{u}(x).
\end{equation*}%
Hence, substituting $\chi \equiv 1$ into (\ref{Rem-2}), we have
\begin{equation*}
\int_{\mathbb{X}}dE_{u,w_{1}}^{(k)}(x)=\frac{u}{1-u}H_{k}^{(r)}(u\mid w_{1}),
\end{equation*}%
which finally yields
\begin{eqnarray*}
&&\int_{\mathbb{X}}\int_{\mathbb{X}}...\int_{\mathbb{X}%
}e^{(x_{1}w_{1}+...+x_{r}w_{r}+w)t}dE_{u}(x_{1})dE_{u}(x_{2})...dE_{u}(x_{r})
\\
&=&\frac{u^{r}}{(e^{w_{1}t}-u)(e^{w_{2}t}-u)...(e^{w_{r}t}-u)}e^{wt}.
\end{eqnarray*}

\begin{acknowledgement}
The present investigation was supported, in part, by the Natural Science and
Engineering Research Council of Canada under Grant OGP0007353.
\end{acknowledgement}


\begin{thebibliography}{99}
\bibitem{Andrews} G. E. Andrews, $q$-analogues of the binomials coefficient
congruences of Babbage, Wolstenhome and Glaiser, \textit{Discrete Math.},
\textbf{204} (1999), 15-25.

\bibitem{Apostol} T. M. Apostol, \textit{Introduction to analytic number
theory}, Springer-Verlag, New York, 1976.

\bibitem{Askey} R. Askey, The $q$-gamma and $q$-beta functions, \textit{%
Appl. Anal}., \textbf{8} (1978), 125-141.

\bibitem{Barnes} W. Barnes, On theory of the multiple gamma functions,
\textit{Trans. Camb. Philos. Soc.}, \textbf{19} (1904), 374-425.

\bibitem{Carlitz} L. Carlitz, $q$-Bernoulli numbers and polynomials, \textit{%
Duke Math.}, \textbf{15} (1948), 987-1000.

\bibitem{Cherednik} I. Cherednik, On $q$-analogues of the Riemann's zeta
function, \textit{Selecta Math.}, \textbf{7} (2001), 447-491.

\bibitem{Deeba} E. Deeba and D. Rodriguez, Stirling's series and Bernoulli
functions, \textit{Amer. Math. Monthly}, \textbf{98} (1991), 423-426.

\bibitem{Dilcher} K. Dilcher, Some $q$-series identities related to divisor
functions, \textit{Discrete Math}., \textbf{145} (1995), 83-93.

\bibitem{E. Friedman and S. Ruijsenaars} E. Friedman and S. Ruijsenaars,
Shintani-Barnes zeta and gamma functions, \textit{Adv. in Math.}, \textbf{187%
} (2004), 362-395.

\bibitem{Howard} F. T. Howard, Applications of a recurrence for the
Bernoulli Numbers, J. Number Theo., 52 (1995), 157-172.

\bibitem{Iwasawa} K. Iwasawa, Lecture on $p$-adic $L$-functions, \textit{%
Annals of Mathematics Studies}, No. \textbf{74}, Princeton University Press,
Princeton NJ, 1972.

\bibitem{M. Jimbo and T. Miwa} M. Jimbo and T. Miwa, Quantum $KZ$ equation
with $\mid q\mid =1$ and correlation functions of the $XXZ$ model in the
gapless regime, \textit{J. Phys. A: Math.}, \textbf{29} (1996), 2923-2958.

\bibitem{Katriel} J. Katriel, Stirling numbers identities interconsistency
of $q$-analogues, \textit{J. Phys. A: Math.}, \textbf{31} (1988), 3559-3572.

\bibitem{Kim1} T. Kim, On explicit formulas of $p$-adic $q$-$L$-functions,%
\textit{\ Kyushu J. Math.}, \textbf{48} (1994), 73-86.

\bibitem{Kim2} T. Kim, On a $q$-analogue of the $p$-adic Log gamma functions
and related integrals, \textit{J. Number Theo.}, \textbf{76} (1999), 320-329.

\bibitem{Kim3} T. Kim, A note on $p$-adic $q$-Dedekind sums, \textit{C. R.
Acad. Bulgare Sc.}, \textbf{54} (2001), 37-42.

\bibitem{Kim4} T. Kim, An invariant $p$-adic integral associated with Daehee
Numbers, \textit{Integral Transform. Spec. Funct.}, \textbf{13} (2002),
65-69.

\bibitem{Kim5} T. Kim, $q$-Volkenborn integration, \textit{Russ. J. Math
Phys.}, \textbf{19} (2002), 288-299.

\bibitem{Kim6} T. Kim, On $p$-adic $q$-$L$-functions and sums of powers,
\textit{\ Discrete Math.}, \textbf{252} (2002), 179-187.

\bibitem{Kim7} T. Kim, Non-archimedean $q$-integrals associated with
multiple Changhee $q$-Bernoulli Polynomials, \textit{Russ. J. Math Phys.},
\textbf{10} (2003), 91-98.

\bibitem{Kim8} T. Kim, $q$-Riemann zeta function, \textit{Internat. J. Math.
Sci.}, \textbf{2004} (2003), 185--192.

\bibitem{Kim9} T. Kim, On Euller-Barnes multiple zeta functions, \textit{%
Russ. J. Math Phys}., \textbf{10} (2003), 261-267.

\bibitem{Kim10} T. Kim, A note on Dirichlet series, \textit{Proc. Jangjeon
Math. Soc}., \textbf{6} (2003),161-166.

\bibitem{Kim11} T. Kim, Sums of powers of consecutive $q$-integers, \textit{%
Adv. Stud. Contep. Math.}, \textbf{9} (2004), 15-18.

\bibitem{Kim12} T. Kim, A note on $q$-zeta functions, \textit{Proceedings of
The 15th International Conference of The Jangjeon Mathematical Society}, (
Hapcheon, South Korea; August 5-7, 2004), 110-114.

\bibitem{Kim13} T. Kim, A note on the $q$-multiple zeta function, \textit{%
Adv. Stud. Contep. Math.}, \textbf{8} (2004), 111-113.

\bibitem{Kim14} T. Kim, $p$-adic $q$-integrals associated with the
Changhee-Barnes' $q$-Bernoulli Polynomials, \textit{Integral Transform.
Spec. Funct.}, \textbf{15} (2004), 415-420.

\bibitem{Kim15} T. Kim, Analytic continuation of multiple $q$-zeta functions
and their values at negative integers, \textit{Russ. J. Math Phys.}, \textbf{%
11} (2004), 71-76.

\bibitem{Kim-Rim} T. Kim and S. -H. Rim, A note on the $q$-integrals and $q$%
-series, \textit{Adv. Stud. Contep. Math.}, \textbf{2} (2000), 37-45.

\bibitem{kim-Jang-Rim-Son} T. Kim, L. C. Jang, S. H. Rim and J. -W. Son, On
the values of zeta and $L$-functions, \textit{Proc. Jangjeon Math. Soc.},
\textbf{1} (2000), 11-18.

\bibitem{T. Kim-SD.Kim-DW.Park} T. Kim, S. D. Kim and D. -W. Park, On
Uniform Differentiability and q-Mahler expansions, \textit{Adv. Stud.
Contep. Math.}, \textbf{4} (2001), 35-41.

\bibitem{T.Kim-S.-H.Rim} T. Kim and S. -H. Rim, On Changhee-Barnes' $q$%
-Euler numbers and polynomials, \textit{Adv. Stud. Contep. Math.}, \textbf{9}
(2004), 81-86.

\bibitem{Khrennikov} A. Khrennikov, $p$\textit{-Adic Valued Distributions in
Mathematical Physics}, Kluwer Academic Publishers, Dordrecht, Boston and
London 1994.

\bibitem{Khrennikov-1} A. Khrennikov, $\ p$-Adic Discrete Dynamical Systems
and Their Applications in Physics and Cognitive Science, \textit{Russ. J.
Math. Phys}., \textbf{11} (2004), 45-70.

\bibitem{N. Koblitz} N. Koblitz, $q$-extension of the $p$-adic gamma
function, \textit{Trans. Amer. Math. Soc}., \textbf{260} (1980), 449-457.

\bibitem{N. Koblitz-1} N. Koblitz, On Carlitz's $q$-Bernoulli numbers,
\textit{J. Number Theory}, \textbf{14} (1982), 332-339.

\bibitem{T. H. Koornwinder} T. H. Koornwinder, Special functions and $q$%
-commuting valuables, \textit{Fields Inst. Comm}., \textbf{14} (1997).

\bibitem{N. Kurokawa} N. Kurokawa, Multiple sine functions and selberg zeta
functions, Proc. \textit{Japan Acad. A}, \textbf{67} (1991), 61-64.

\bibitem{K. Matsumoto} K. Matsumoto, The analytic continuation and the
asyptotic behaviour of certain multiple zeta-function I, \textit{J. Number
Theory,} \textbf{101} (2003), 223-243.

\bibitem{C. A. Nelson and M. G. Gartley} C. A. Nelson and M. G. Gartley, On
the zeros of the $q$-analogue exponential function, \textit{J. Phys. A: Gen.
Math}., \textbf{27} (1994), 3857-3881.

\bibitem{C. A. Nelson and M. G. Gartley-1} C. A. Nelson and M. G. Gartley,
On the two $q$-analogue logaritmic functions: $\ln _{q}(w),\ln (\ln _{q}(z))$%
, \textit{J. Phys. A: Gen. Math.}, \textbf{24} (1996), 8099-8115.

\bibitem{M. Nishizawa} M. Nishizawa, On a $q$-analogue of the multiple gamma
functions, \textit{Lett. Math. Phys.}, \textbf{37} (1996), 2001-2009.

\bibitem{K. Ota} K. Ota, On Kummer-type congruences for derivatives of
Barnes' multiple Bernoulli Polynomials,\textit{\ J. Number Theory,}\textbf{\
92 }(2002), 1-36.

\bibitem{K. Ota-1} K. Ota, Derivatives of Dedekind sums and their
reciprocity law, \textit{J. Number Theory}, \textbf{98} (2003), 280-309.

\bibitem{T. M. Rassia and H. M. Srivastava} T. M. Rassia and H. M.
Srivastava, Some classes of infinite series associated with the Riemann zeta
and polygamma functions and generalized harmonic numbers, \textit{Appl.
Math. Computation}, \textbf{131} (2002), 593-605.

\bibitem{A. M. Robert} A. M. Robert, \textit{A course in }$p$\textit{-adic
Analysis}, Springer-Verlag, New York, 2000.

\bibitem{S. N. M. Ruijsenaars} S. N. M. Ruijsenaars, On Barnes' multiple
zeta function and gamma functions, \textit{Adv. in Math.}, \textbf{156}
(2000), 107-132.

\bibitem{B. E. Sagan} B. E. Sagan, Congruence properties of $q$analogues,
\textit{Adv. in Math.}, \textbf{95} (1992), 127-143.

\bibitem{W. H. Schikhof} W. H. Schikhof, \textit{Ultrametric Calculus},
Cambridge University Press, Cambridge, London and New York, 1984.

\bibitem{K. Shiratani} K. Shiratani, On Euler numbers, \textit{Mem. Fac.
Kyushu Uni}., \textbf{27} (1973), 1-5.

\bibitem{K. Shiratani and S. Yamamoto} K. Shiratani and S. Yamamoto, On a $p$%
-adic interpolation function for the Euler numbers and its derivative,
\textit{Mem. Fac. Kyushu Uni.}, \textbf{39} (1985), 113-125.

\bibitem{Y. Simsek} Y. Simsek, Theorems on twisted $L$-functions and twisted
Bernoulli numbers, to appear \textit{Proc. Jangjeon Math. Soc.}

\bibitem{Y. Simsek-1} Y. Simsek, Generalized Dedekind sums associated with
the Abel sum and the Eisenstein and Lambert series, \textit{Adv. Stud.
Contep. Math.} \textbf{9} (2) (2004), 125-137.

\bibitem{H. M. Srivastava and J. Choi} H. M. Srivastava and J. Choi, \textit{%
Series Associated with the Zeta and Related Functions}, Kluwer Acedemic
Publishers, Dordrecht, Boston and London, 2001.

\bibitem{H. M. Srivastava and P. W. Karlsson} H. M. Srivastava and P. W.
Karlsson, \textit{Multiple Gaussian Hypergeometric Series}, Halsted Press
(Ellis Horwood Limited, Chichester), John Wily and Sons, New York,
Chichester, Brisbane and Toronto, 1985.

\bibitem{H. Tsumura-1} H. Tsumura, On a $p$-adic interpolation of
generalized Euler Numbers and its applications, \textit{Tokyo J. Math}.,
\textbf{10} (1987), 281-293.

\bibitem{P. T. Young} P. T. Young, On the behavior of some two-variable $p$%
-adic $L$-function, \textit{J. Number Theory}, \textbf{98} (2003), 67-86.

\bibitem{Vilademir} V. S. Vilademir, I. V. Volvoich and E. I. Zelenov,%
\textit{\ }$p$\textit{-adic Analysis and Mathematical Physics}, Cambridge
Univ. Press, London and New York,1990.

\bibitem{Washington} L. C. Washington, \textit{Introduction to Cyclomotic
Fields}, Springer-Verlag and New York, 1997.

\bibitem{E. T. Wittaker and G. N. Watson} E. T. Wittaker and G. N. Watson,
\textit{A course of modern Analysis}, Cambridge Univ. Press, London and New
York, 1927.

\bibitem{S. C. Woon} S. C. Woon, Fractal of the Julya and Mandelbort sets of
the Riemann zeta functions, (arXiv: Chao-dyn 19812031VI 27 Dec. 1998).
Preprint(1998).
\end{thebibliography}
\end{document}